\documentclass[reqno,11pt,twoside]{amsart}
\usepackage{amsmath}
\usepackage{amsthm}
\usepackage{graphicx}
\usepackage{amsfonts}
\usepackage{amssymb}
\usepackage{bm}
\usepackage{color}
\usepackage{hyperref}
\usepackage{enumerate}
\usepackage[margin=1.2in]{geometry}
\usepackage{color,soul}
\usepackage{hyperref}

\numberwithin{equation}{section}

\newtheorem{theorem}{Theorem}[section]
\newtheorem{lemma}[theorem]{Lemma}

\newtheorem{proposition}[theorem]{Proposition}

\theoremstyle{definition}

\newtheorem{definition}[theorem]{Definition}
\newtheorem{assumption}[theorem]{Assumption}
\newtheorem{remark}[theorem]{Remark}

\def\Beta{\mathfrak{B}}
\def\E{{\mathbb E}}

\def\R{{\mathbb R}}

\def\N{{\mathbb N}}
\def\PP{{\mathbb P}}

\def\P{{\mathcal P}}

\def\L{{\mathcal L}}

\def\W{{\mathcal W}}

\def\F{{\mathcal F}}
\def\FF{{\mathbb F}}

\def\dd{\mathrm{d}}

\DeclareMathOperator*{\sgn}{sgn}

\title{From Rank-based models with common noise to Pathwise entropy solutions of SPDEs}
\author{Mykhaylo Shkolnikov and Lane Chun Yeung}
\address{Department of Mathematical Sciences, Carnegie Mellon University}
\email{mshkolni@gmail.com, lyeung@andrew.cmu.edu}
\thanks{M. Shkolnikov is partially supported by the NSF grant DMS-2108680.}

\usepackage{fancyhdr}
\begin{document}

\begin{abstract}
	We study the mean field limit of a rank-based model with common noise, which arises as an extension to models for the market capitalization of firms  in stochastic portfolio theory. We show that, under certain conditions on the drift and diffusion coefficients, the empirical cumulative distribution function converges to the solution of a stochastic PDE. A key step in the proof, which is of independent interest,  is to show that any solution to an associated martingale problem is also a pathwise entropy solution to the stochastic PDE, a notion introduced in a recent series of papers \cite{lions2013scalar,lions2014scalar,hofmanova2016scalar,gess2017long,gess2017stochastic}.  
\end{abstract}
\maketitle

\section{Introduction}

We study the following system of interacting diffusion processes on the real line:
\begin{align}\label{eq: general rank based system}
\begin{split}
\dd X^{n,i}_t &=  b\big(F_{\nu^n_t}\big(X^{n,i}_t\big)\big) \,  \dd t +  \sigma\big(F_{\nu^n_t}\big(X^{n,i}_t\big)\big) \, \dd B^{n,i}_t + 
\gamma\big(F_{\nu^n_t}\big(X^{n,i}_t\big)\big) \, \dd W^{n}_t, 
\end{split}
\end{align} 
for $ i = 1, \dots, n $ and $ 0 \le t \le T$. Here, $ \nu^n_t := \frac{1}{n} \sum_{i=1}^n \delta_{X^{n,i}_t} $  denotes  the empirical measure of the particle system at time $ t $, $ F_{\nu^n_t} $ is its cumulative distribution function (CDF),  $ b : [0,1] \to \R$ and  $ \sigma, \gamma : [0,1] \to (0, \infty)$ are measurable functions, and $B^{n,1}, \dots, B^{n,n}, W^{n}$ are independent one-dimensional standard Brownian motions. The system \eqref{eq: general rank based system} is a \emph{rank-based}  model because the drift and diffusion coefficients of each particle are determined by its rank in the particle system.

Models with piecewise constant coefficients arose originally from questions in filtering theory \cite{bass1987uniqueness}. More recently, there has been a lot of interest in  rank-based models without common noise, i.e., when $\gamma \equiv 0$,  due to their applications in stochastic portfolio theory (see, e.g.,  \cite{fernholz2002stochastic, banner2005atlas,  karatzas2009stochastic,ichiba2013convergence,almada2019dynamics,itkin2022robust,itkin2021class,itkin2024calibrated}).  In this context, $X^{n,i}$ represents the logarithmic market capitalization of the $i$-th firm, and models of this form are known to be able to capture some structures of real financial markets, in particular the shape and stability of the capital distribution curve \cite{fernholz2002stochastic,banner2005atlas,chatterjee2010phase, pal2014concentration,pal2008one}.
 This observation naturally leads  to  the study of the large $n$ behavior of such models \cite{shkolnikov2012large,dembo2016large,kolli2018spde}, which describes the dynamics of the capital distribution among a large number of companies, and is therefore of particular interest to institutional investors whose portfolios typically include stocks of  hundreds or thousands of companies.

One shortcoming of rank-based models without common noise is that  only idiosyncratic noises drive the stock prices. This means that a firm's stock is only exposed to idiosyncratic risk, i.e., an inherent risk that affects only that specific firm, such as poor sales of a particular product, or change in management. A richer model would also allow for  systematic risk that affects the entire market, such as change in interest rates, inflation, or other macroeconomic factors. To this end,  \cite{kolli2019large} took a first step and incorporated an additional term $\gamma(t, \nu^n_t) \, \dd W^{n}_t$ common to all stocks in the model. For that model, they proved in \cite[Theorem 1.2]{kolli2019large} a version of the law of large numbers (``hydrodynamic limit") asserting that as $n \to \infty$, the empirical CDF converges to a limit which is the solution of a stochastic PDE. It is worthwhile to point out that their common noise term affects all particles equally. On the other hand, while the common noise term in \eqref{eq: general rank based system} does not have a time dependence, it allows for particles to experience different effects from the common noise depending on their ranks. In some sense, this is a more realistic model because a firm's size is known to be negatively correlated with its beta coefficient in the Capital Asset Pricing Model \cite{sullivan1978cost,subrahmanyam1980systematic,binder1992beta}. However, this new form of rank-based models brings significant mathematical challenges when it comes to proving the corresponding law of large numbers, as it necessitates  a study of the relationship between the pathwise entropy solution of the limiting stochastic PDE and the solution of an associated martingale problem (see Subsection \ref{subsec:mp.pe} for a detailed discussion). More recently,  \cite{cuchiero2024polynomial} proved the corresponding law of large numbers for volatility-stabilized market models, which is  another class of models in stochastic portfolio theory.  

From another point of view, the system \eqref{eq: general rank based system} without common noise  falls under the general framework of  mean field interacting diffusions, which originates from the seminal work of McKean \cite{mckean1967propagation}. In this regard, the drift and diffusion coefficients of each particle are viewed as functions of the  current position of the particle $ X^{n,i}_t $ and the empirical measure of the  system $ \nu^n_t $. Assuming that the drift and diffusion coefficients are jointly continuous (in the state and measure variables),  G\"{a}rtner established in \cite{gartner1988mckean} a law of large numbers (which is also called propagation of chaos or convergence to the mean field limit in the literature). Unfortunately, even without common noise, his result is not applicable to \eqref{eq: general rank based system} because the drift and diffusion coefficients are discontinuous in both the position of the particle and the empirical measure of the system. Nonetheless, the special structure of the coefficients permits a derivation of the law of large numbers, see \cite[Proposition 2.1]{jourdain2013propagation}, \cite[Theorem 1.4]{jourdain2000probabilistic}, \cite[Theorem 1.2]{shkolnikov2012large} and \cite[Corollary 1.6]{dembo2016large}, where it is shown that the empirical CDF converges to the solution of the porous medium PDE.  Moreover, a  central limit theorem concerning the fluctuations of the empirical CDF around its limit was proven in \cite[Theorem 1.2]{kolli2018spde}, and a large deviations principle was obtained in \cite[Theorem 1.4]{dembo2016large}.

In this paper, we show that, under suitable assumptions on $ b $, $ \sigma $ and $ \gamma $, and on the  initial positions of the particles, as $ n \to \infty $, the empirical measure process $\nu^n$ of \eqref{eq: general rank based system} converges in distribution to a limit $\nu$, whose CDF  process $u(t, \cdot) := F_{\nu_t}(\cdot)$ solves the stochastic PDE
\begin{align} \label{SPDE ito}
	\dd u = \left(-\Beta(u)_x + \Sigma(u)_{xx} + \Gamma(u)_{xx}\right) \dd t - G(u)_x \, \dd W_t, \quad t \in [0,T],
\end{align}
where $W$ is a one-dimensional standard Brownian motion,  and $\Beta : [0,1] \to \R$ and $\Sigma, \Gamma, G : [0,1] \to [0, \infty)$ are defined by
\begin{align} \label{eq: def B.Sigma.Gamma}
	\begin{split}
		\Beta(r) &:= \int_0^r b(a) \, \dd a, \qquad \Sigma(r) := \frac{1}{2}\int_0^r \sigma^2(a) \, \dd a,\\
		\Gamma(r)&:= \frac{1}{2} \int_0^r \gamma^2(a) \, \dd a, \qquad   G(r) := \int_0^r \gamma(a) \, \dd a.
	\end{split}
\end{align}
 Note that the SPDE \eqref{SPDE ito}  describes the evolution of the conditional CDF in the McKean-Vlasov equation
\begin{align*}
\dd X_t = b\left(F_{\nu_t}(X_t)\right) \dd t+ \sigma\left(F_{\nu_t}(X_t)\right) \dd B_t + \gamma\left(F_{\nu_t}(X_t)\right) \dd W_t, \quad \nu_t = \L(X_t | W_{[0,t]}),
\end{align*}
which intuitively is the large $n$ limit of the particle system \eqref{eq: general rank based system}. Here, $ \L(X_t | W_{[0,t]})$ denotes (a version of) the conditional law of $ X_t $ given the path of the Brownian motion $ W $ on $ [0, t] $. 
\subsection{Notation}

We employ the usual notation $ \langle \nu, f \rangle = \int_\R f \, \dd \nu$ to denote  the integration of a real-valued $\nu$-integrable function $f$ with respect to a Borel measure $\nu$ on $ \R $. Similarly, for a vector of $\nu$-integrable functions $\boldsymbol{f} = (f_1, \dots, f_k) $, we write $\langle \nu, \boldsymbol{f} \rangle = (\langle \nu, f_1\rangle, \dots, \langle \nu, f_k\rangle) \in \R^k$. For two real-valued functions $ g $ and $ f $, we use the notation $ \langle g,f \rangle = \int_\R g\, f \, \dd x $ whenever $g f \in L^1(\R)$, and write $\langle g, \boldsymbol{f} \rangle = (\langle g, f_1\rangle, \dots, \langle g, f_k\rangle) $ for any vector of functions  $\boldsymbol{f} = (f_1, \dots, f_k) $. For a real-valued function $f$ on any set $E$, we let $\Vert f \Vert_\infty := \sup_{x \in E} |f(x)|$.  

For a metric space $(E,d)$, let $\P(E)$ denote the space of Borel probability measures on $E$, equipped with the topology of weak convergence. Let $\P_1(\R)$ denote the subspace of $ \P(\R) $ with finite first moments, equipped with the Wasserstein distance $\W_1$ defined by
\begin{align}\label{eq:def.W1}
	\W_1(\mu, \nu) = \inf_{\pi} \int_{\R \times \R} |x - y|\, \pi(\dd x, \dd y).
\end{align}
Here, the infimum is taken over all $\pi \in \P(\R \times \R)$ with marginals  $\mu$ and $\nu$.  

Let $ T > 0 $ be fixed throughout the paper. For $0 \le s \le t\le T$ and a metric space $(E,d)$, we use $ C([s,t]; E)$ to denote the space of continuous functions $x: [s,t] \to E$ equipped with the topology of uniform convergence. When $E = \R$, this topology is the same as the one induced by the uniform norm
\begin{align*}
	\Vert x \Vert_{C([s,t];\R)} := \sup_{r \in [s,t]} |x_r|.
\end{align*}  
We also set $C_0([0,T];\R)  : = \{x \in C([0,T]; \R) : x(0) = 0\}$. 

The subspace of $L^\infty(\R)$-functions with bounded total variation is denoted by $BV(\R)$, and we use $\Vert f \Vert_{BV(\R)}$ to denote the total variation of $f \in BV(\R)$. Finally, for any $F:[0,T] \to \R$, we let $F\big|_s^t := F(t) - F(s)$. 
\subsection{Main result}
We make the following two assumptions.
\begin{assumption} \label{assum}
		\begin{enumerate}[(a)]
		\item $ b :[0,1] \to \R $ is continuously differentiable.  
		\item $ \gamma :[0,1] \to (0, \infty)$ is continuously differentiable.   \label{eq: assum gamma}
		\item $\sigma : [0, 1] \to (0, \infty)$ is bounded and  continuous.  Also, it is non-degenerate: $\inf_{a \in [0,1]} \sigma(a) > 0$. \label{eq: assum sigma}
	\end{enumerate}
\end{assumption}
\begin{assumption} \label{assum:init}
There exists a $ \nu^0 \in \P_1(\R) $ such that $ \W_1(\nu^n_0, \nu^0) \to 0$ in distribution.
\end{assumption}
\begin{remark}
	Assumption \ref{assum:init} is satisfied if the initial positions $X^{n,1}_0, \dots, X^{n,n}_0$ are i.i.d. with some common distribution $\nu^0 \in \P_1(\R)$. Indeed, $\nu^n_0 \to \nu^0 $ weakly a.s.  by Varadarajan's Theorem  \cite[Theorem 11.4.1]{dudley2018real},  and $\frac{1}{n} \sum_{i=1}^n |X^{n,i}_0| \to \int_\R |x| \, \nu^0(\dd x)$ a.s. by the law of large numbers. Therefore,  by \cite[Theorem 7.12 (iii) $\Rightarrow$ (i)]{villani2003topics}, we deduce $  \W_1(\nu^n_0, \nu^0) \to 0$ a.s., which implies $  \W_1(\nu^n_0, \nu^0) \to 0$ in distribution. 
\end{remark}
Our main result can now be stated as follows. 
\begin{theorem} \label{thm: rank}
	Suppose Assumptions \ref{assum} and \ref{assum:init} hold. 
Then for each $n \in \N$, there exists a weak  solution to \eqref{eq: general rank based system}, which is unique in law. 
	The sequence $\left(\nu^n\right)_{n \in \N}$ of  $C([0,T];\P_1(\R))$-valued random variables converges in distribution to a $C([0,T];\P_1(\R))$-valued random variable $\nu$, such that the corresponding CDFs $u := (u(t,\cdot))_{t \in [0,T]} := (F_{\nu_t})_{t \in [0,T]}$ form a weak solution (in both the probabilistic and PDE sense) to the SPDE \eqref{SPDE ito}, i.e., there exists a  one-dimensional standard Brownian motion $ W $ such that for all  $0 \le s \le t \le T$ and $f \in C^\infty_c(\R)$, 
	\begin{equation} \label{eq: weak sol SPDE}
		\begin{split}
			&\int_{\R}  u(t,x) \, f(x) \, \dd x -   \int_{\R}  u(s,x) \, f(x) \,\dd x  =  \int_s^t \int_{\R} G(u(r,x)) \, f'(x) \, \dd x \, \dd W_r\\
			 &  \hspace{0.5cm}+ \int_s^t \int_{\R} \Beta(u(r,x)) \, f'(x) +   \Sigma(u(r,x))  \, f''(x) + \Gamma(u(r,x)) \, f''(x) \,  \dd x \, \dd r, 
		\end{split}
	\end{equation}
	a.s., with initial condition $ u(0,\cdot) =  F_{\nu^0}(\cdot) $. Moreover,  pathwise uniqueness holds for the SPDE \eqref{SPDE ito}. In particular, the law of $ u $ is unique. 
\end{theorem}
\begin{remark}
	With a bit more bookkeeping, Theorem \ref{thm: rank}  extends to the case of multiple common noises, i.e., to models of the form 
	\begin{align*}
		\begin{split}
			\dd X^{n,i}_t &=  b\big(F_{\nu^n_t}\big(X^{n,i}_t\big)\big) \, \dd t+  \sigma\big(F_{\nu^n_t}\big(X^{n,i}_t\big)\big) \, \dd B^{n,i}_t + 
			\sum_{j=1}^{k} \gamma_j\big(F_{\nu^n_t}\big(X^{n,i}_t\big)\big) \, \dd W^{n,j}_t,
		\end{split}
	\end{align*}  
	where $W^{n,j}$ are independent standard Brownian motions and $\gamma_j : [0,1] \to (0, \infty)$ are continuously differentiable. Since the extension is straightforward, we  focus on the case with just one common noise in this paper. 
\end{remark}
\subsection{Martingale problem and pathwise entropy solution} \label{subsec:mp.pe}
The proof of Theorem \ref{thm: rank} follows the well-trodden path of tightness-limit-uniqueness.  Central to the uniqueness proof is that any solution to a martingale problem associated with \eqref{SPDE ito}, in the sense of Stroock and Varadhan \cite{stroock1997multidimensional},  is also a \emph{pathwise entropy solution}. 
To explain the latter notion of solution, let us first recast \eqref{SPDE ito} into a  more general form:
\begin{align} \label{eq: multi SPDE}
\begin{cases}
\dd u = \big(-\Beta(u)_{x} + \Sigma(u)_{xx} \big) \, \dd t -  G(u)_{x} \circ \dd z, \quad \text{on} \quad [0, T] \times \R,\\
u= u^0,   \quad \text{on} \quad  \{0\} \times  \R  ,	
\end{cases}
\end{align}
where $ u^0 \! \!: \!\R \to \R$, $ z \in C_0([0,T]; \R)$,   $\Beta$, $\Sigma$, and $G$ are as given in \eqref{eq: def B.Sigma.Gamma}, and $\circ \, \dd z$ denotes the Stratonovich differential. When the driving signal $z$ is the one-dimensional standard Brownian motion,  \eqref{eq: multi SPDE} is the Stratonovich formulation of \eqref{SPDE ito}.  

The multidimensional and driftless version (i.e., when  $ \Beta \equiv 0 $)  of this stochastic degenerate parabolic-hyperbolic equation  is  studied in \cite{gess2017stochastic}, building on earlier developments in \cite{lions2013scalar,lions2014scalar,hofmanova2016scalar,gess2017long}. There, the notion of \emph{pathwise entropy solution} for \eqref{eq: multi SPDE} is introduced, which is based on evaluating test functions for the ``kinetic formulation" of \eqref{eq: multi SPDE}  along the characteristics of a suitable transport equation.   
We outline the key ideas in the construction of this notion of solution, and refer the reader to \cite{gess2017stochastic}  for more details. The construction starts from the kinetic formulation of \eqref{eq: multi SPDE}: Define    $\overline{\chi} \!:  \R \times \R \to \R $ by
\begin{align} \label{chi}
	\overline{\chi} (\xi, u ) = \begin{cases}
		1, \quad \text{if  } 0 < \xi < u,\\
		-1, \quad \text{if  } u < \xi < 0,\\
		0, \quad \text{otherwise},
	\end{cases}
\end{align}
and, given $u: [0, T] \times \R \to \R$,  let 
\begin{align*}
	\chi(\xi, t, x) :=\overline{\chi} \big(\xi, u(t,x) \big). 
\end{align*}
Then the \emph{kinetic formulation} of \eqref{eq: multi SPDE} is the following Cauchy problem for $ \chi $:
\begin{align}  \label{kinetic}
	\begin{cases}
	\partial_t \chi + \big(b(\xi) + \gamma(\xi) \dot{z}_t\big) \,  \partial_x \chi - \frac{1}{2} \sigma^2(\xi) \, \partial^2_{xx} \chi  = \partial_\xi (m + n) \quad \text{on} \quad \R \times [0, T]  \times \R, \\
	\chi(\xi , 0, x	) = \overline{\chi}(\xi, u^0(x)),
	\end{cases}
\end{align}
where the ``entropy defect measure" $ m $ and the ``parabolic dissipation measure" $ n $ are non-negative finite measures on $ \R \times [0,T] \times \R$, and $ \dot{z}_t $ denotes the time derivative of $ z_t $. 

Observe that for \eqref{kinetic} to make sense, the driving signal $ z $ needs to be differentiable, which rules out the case of Brownian motion. However, the remarkable observation in \cite{gess2017stochastic} is that by carefully choosing  a set of test functions for \eqref{kinetic}, the terms involving $ \dot{z}_t $ can be eliminated. More specifically,  let us consider the  transport equation
\begin{align} \label{eq: transport.eq multi}
		\partial_t \varrho(\xi, t, x) + \left(b(\xi) + \gamma(\xi) \dot{z}_t\right)   \varrho_{x}(\xi, t, x)  = 0, 
		\quad \text{on} \quad \R \times [0,T] \times \R.
\end{align}
For each $ (\eta, y) \in \R^{2}$ and $ \varrho^0 \in C^\infty_c(\R^2)$, note that $ t \mapsto x +  y + b(\xi)t + \gamma(\xi) z_t  $ are   characteristics  of \eqref{eq: transport.eq multi}, and so  
\begin{align}  \label{eq: varrho multi}
	\varrho(\xi, t, x; \eta,  y) = \varrho^0(x -y -b(\xi) t - \gamma(\xi)z_t, \xi - \eta)
\end{align}
 is a solution to \eqref{eq: transport.eq multi}, with the initial condition  $\varrho^0(x-y, \xi - \eta)$. It can then be shown (cf. \cite[Lemma 2.2]{gess2017stochastic}) that when \eqref{kinetic} is tested against functions of the form \eqref{eq: varrho multi}, one obtains
 	\begin{equation} 
 	\label{eq: pathwise entropy solution result multi}
 	\begin{split}
 		&- \int_{\R^{2}} \overline{\chi}(\xi, u(\cdot ,x))  \, \varrho(\xi, \cdot, x; \eta, y) \, \dd \xi \, \dd x \bigg|_s^t+  \frac{1}{2} \int_s^t \int_{\R^2} \overline{\chi} \left(\xi, u(r,x)\right) \,  \sigma^2(\xi) \, \varrho_{xx}(\xi, r, x ; \eta ,y) \, \dd \xi \, \dd x \, \dd r\\
 		&=   \int_s^t \int_{\R^2} \partial_{\xi}  \varrho(\xi, r, x ; \eta ,y) \, (m + n)(\dd x , \dd \xi, \dd r).
 	\end{split}
 \end{equation} 
 In particular, even though the smoothness of $z$ was used to derive \eqref{eq: pathwise entropy solution result multi}, the identity   \eqref{eq: pathwise entropy solution result multi} itself makes sense for all  $z \in C_0([0,T]; \R)$. This identity forms the basis of the definition of pathwise entropy solutions, which we  give next. In order to tailor it to our setting, the definition is slightly modified from \cite[Definition 2.1]{gess2017stochastic}. See Remark \ref{rmk:def.pe} for further discussion. We let
\begin{align} \label{eq: beta multi}
S(r) :=  \int_0^r \sigma(a)\, \dd a.
\end{align}
\begin{definition} \label{def:pe}
	Let $ u^0 \in L^\infty(\R) $. A function $ u \in  L^\infty([0,T] \times \R) $ satisfying 
	\begin{align}  \label{u.continuity.0}
		 \int_\R |u(s,x) - u(t,x) | \, \dd x  < \infty, \quad \text{for all } 0 \le s \le t \le T
	\end{align}
	and 
	\begin{align} \label{u.continuity}
		\lim_{s \to t} \int_\R |u(s,x) - u(t,x) | \, \dd x = 0, \quad \text{for all } t \in [0,T] 
	\end{align}
	 is a \emph{pathwise entropy solution} to  \eqref{eq: multi SPDE} with respect to $ z $ and with initial condition $ u^0 $ if $ u(0, \cdot) = u^0 $ and 
	\begin{enumerate}[(a)]
		\item \begin{align} \label{PE.def.1}
			  S(u)_x \in L^2([0,T] \times \R ). 
		\end{align}
		\item  For all test functions $\varrho $ given by \eqref{eq: varrho multi} with $\varrho^0 \in C^\infty_c(\R^2)$ and all $ (\eta, y) \in \R^2$,
		\begin{align} \label{eq: chain rule}
		 \int_{\R^{2}} \overline{\chi}(\xi, u(t,x))  \, \sigma(\xi) \, \varrho_x(\xi, t, x; \eta, y) \,\dd x\, \dd \xi= -  \int_{\R} S \left(u(t,x)\right)_{x} \varrho(u(t,x), t, x; \eta ,y) \, \dd x
		\end{align}
		holds for a.e. $ t \in [0, T] $.
		\item  There exists a non-negative finite measure $m$ on $  \R^2 \times [0,T] $ such that for all test functions $\varrho $ given by \eqref{eq: varrho multi} with $\varrho^0 \in C^\infty_c(\R^2)$,  all $ (\eta, y) \in \R^2$ and all $0 \le s \le t \le T$,
the identity \eqref{eq: pathwise entropy solution result multi} holds 
with $n$ being the non-negative finite measure on $ \R^2 \times [0,T] $ defined by
	\begin{align} \label{eq: measure n}
			n( \dd x,\dd \xi, \dd t ) &:=   \frac{1}{2}\left(  S\left(u(t,x)\right)_x\right)^2   \,  \delta_{u(t,x)}  (\dd \xi) \, \dd x\,  \dd t.
	\end{align}
\end{enumerate} 
\end{definition}
\begin{remark} \label{rmk:def.pe}
	The above definition of pathwise entropy solutions differs from   \cite[Definition 2.1]{gess2017stochastic} in two ways. Firstly, the $ L^1 $-integrability conditions on $ u_0 $ and $ u $ are removed. This is necessary to accommodate the case that $ u(t, \cdot) $ is a CDF, which is the focus of this paper. Because of this, the original continuity requirement $u \in C([0, T]; L^1(\R))$ is changed to \eqref{u.continuity.0} and \eqref{u.continuity}. Secondly, the ``chain rule" \eqref{eq: chain rule} is required to hold only for  a.e. $ t \in [0,T] $ instead of for all $ t \in [0,T] $. This is a minor change in order to accommodate the case that $ u(t,\cdot) $ is continuous only at Lebesgue a.e. $ t \in [0,T] $. 
\end{remark}

We have the following uniqueness and stability result for pathwise entropy solutions in our setting.  It can be seen as an extension of \cite[Theorem 2.3]{gess2017stochastic} in the one-dimensional case. 
\begin{proposition} \label{prop:GSextension}
	Assume that $\sigma $ is positive and bounded, $b$ is continuously differentiable and $\gamma $ is positive and continuously differentiable.  Let $u^{(1)}, u^{(2)} \in L^\infty([0,T] ; BV(\R))$ be two pathwise entropy solutions to \eqref{eq: multi SPDE} with driving signals $z^{(1)}, z^{(2)} \in C_0([0,T]; \R)$ and initial values $ u^{1}_0, u^{2}_0 \in  BV(\R)$. Let $m^{(2)}$ and $n^{(2)}$ denote the finite measures on $\R^2 \times [0, T]$ corresponding to $u^{(2)}$ as given in Definition \ref{def:pe}(c), and $q^{(2)} = m^{(2)} + n^{(2)}$. Then for all $ 0 \le s \le t \le T $, there exists a $ C <\infty  $, which may depend on $\Vert u^{(1)} \Vert_{L^\infty([s,t] ;  BV(\R))}$, $\Vert u^{(2)} \Vert_{L^\infty([s,t] ; BV(\R))}$ and $q^{(1)}(\R^2 \times [s,t])$,  $q^{(2)}(\R^2 \times [s,t])$, such that
	\begin{align*}
	\Vert u^{(1)}(t, \cdot) - u^{(2)}(t, \cdot)\Vert_{L^1(\R)} \le  & \, \Vert u^{(1)}(s,\cdot) - u^{(2)}(s,\cdot)\Vert_{L^1(\R)} + C \Vert z^{(1)} - z^{(2)} \Vert^{1/2}_{C([s,t]; \R)} \\
	&+ C \Vert z^{(1)} - z^{(2)} \Vert_{C([s,t]; \R)} .
	\end{align*}
\end{proposition}
With Proposition \ref{prop:GSextension} in place, the main ingredient in the proof of the uniqueness part of Theorem  \ref{thm: rank}  is the following theorem, which says that  under Assumption \ref{assum}, any  solution to a martingale problem associated with \eqref{SPDE ito} is also a pathwise entropy solution.
 
\begin{theorem} \label{thm: main}
			Suppose Assumption \ref{assum} holds. Let $ (\Omega, \F, \FF= (\F_t)_{t \in [0,T]}, \PP) $ be a filtered probability space and $  \nu \in C([0,T];  \P_1(\R))$ be an  $ \FF $-adapted probability measure-valued process, with $ \nu_0 $ being deterministic. Let $ u(t, \cdot) := F_{\nu_t}(\cdot)$. Suppose $\PP$-a.s., for all $ k \in \N $,  $\boldsymbol{f} = (f_1, \dots, f_k) \in C^\infty_c(\R)^k$, and $\phi \in C^\infty_c(\R^k)$, the process 
	\begin{equation} \label{eq: MRT condition}
		\begin{split}
				[0,T] \ni t \mapsto  & 	\, \phi \left(	\langle u(t,\cdot),  \boldsymbol{f}\rangle \right)  - \phi \left(	\langle  u(0,\cdot), \boldsymbol{f}\rangle \right)  \\
				&- 	\sum_{i=1}^k \int_0^t     \partial_i\phi \left(	\langle  u(r,\cdot), \boldsymbol{f}\rangle \right)  \Big(\big\langle \Beta(u(r,\cdot)), f_i'\big\rangle+ \big\langle (\Sigma + \Gamma)(u(r,\cdot)), f_{i}''\big\rangle \Big) \,  \dd r \\
			\quad & - \frac{1}{2}\sum_{i,j=1}^k \int_0^t  \partial_{ij} \phi\left(	\langle u(r,\cdot), \boldsymbol{f} \rangle \right)  \left\langle G(u(r,\cdot)), f_i' \right\rangle \left\langle G(u(r,\cdot)), f_j' \right\rangle  \, \dd r  
		\end{split}
	\end{equation}
	is an $\FF$-martingale.
	Then:
		\begin{enumerate}[(i)]
			\item \label{mainthm: mp.spde} There exists  an extension $ (\widetilde{\Omega}, \widetilde{\F}, \widetilde{\FF} = (\widetilde{\F}_t)_{t \in [0,T]}, \widetilde{\PP}) $ of the probability space $ (\Omega, \F, \FF = (\F_t)_{t \in [0,T]}, \PP) $ supporting a one-dimensional standard  $\widetilde{\FF}$-Brownian motion $ W $ such that for all  $0 \le s \le t \le T$ and $f \in C^\infty_c(\R)$, \eqref{eq: weak sol SPDE} holds $\widetilde{\PP}$-a.s.
		\item \label{mainthm: ac} $\widetilde{\PP}$-a.s., for a.e. $ t \in [0,T] $, $ \nu_t $ has a density,  or equivalently, $u(t, \cdot)$ is absolutely continuous as a function.
		\item \label{mainthm: u_x in L2}  $ u_x \in L^2([0,T] \times \R)$, $\widetilde{\PP}$-a.s.
		\item \label{mainthm: pe}$\widetilde{\PP}$-a.s.,  $ u $ is a pathwise entropy solution to \eqref{SPDE ito}  with respect to $W$ and with initial condition $ F_{\nu_0}(\cdot)$. 
	\end{enumerate}
\end{theorem}
\subsection{Organization of the paper} The rest of the paper is structured as follows. In Section \ref{sec:proof.mp.pe}, we prove Theorem \ref{thm: main}, relying in particular on a careful study of \eqref{eq: pathwise entropy solution result multi} from a stochastic analysis perspective. In Section \ref{sec:pf.main.thm}, we prove Theorem \ref{thm: rank} by first establishing tightness of the empirical measures, and that any limit point solves the martingale problem described in Theorem \ref{thm: main}.  For the latter, we use techniques similar to \cite[proof of Lemma 1.5]{jourdain2000probabilistic}. Subsequently, we invoke Theorem \ref{thm: main} and Proposition   \ref{prop:GSextension}  to derive the desired uniqueness. The proof of Proposition \ref{prop:GSextension} is given in  Appendix \ref{sec:proof.gs}, where we highlight the differences with the original proof of \cite[Theorem 2.3]{gess2017stochastic}.  Some auxiliary results used in Section \ref{sec:proof.mp.pe} are provided in Appendix \ref{sec:aux}.  
\section{Proof of Theorem \ref{thm: main}} \label{sec:proof.mp.pe}
In this section, we prove Theorem \ref{thm: main}  in several steps. For any $\varepsilon > 0$, set 
\begin{align} \label{eq: varphi}
 \varphi_{\varepsilon}: \R \to \R, \quad 	x \mapsto  \frac{1}{\sqrt{2 \pi \varepsilon}} \exp\left(-\frac{x^2}{2 \varepsilon}\right). 
\end{align}
For any function $ f: [0,T] \times \R \to \R $ and $ t \in [0,T] $, let $ f^\varepsilon (t,\cdot):= f(t, \cdot) \ast \varphi_{\varepsilon}$ denote its convolution with $ \varphi_\varepsilon$ in the spatial variable.  Since convolution commutes with  differentiation, the notation $ \partial_x f^\varepsilon (t,\cdot)$ is unambiguous.   In the following, unless mentioned otherwise, all statements up to \eqref{def:Nd} are to be understood in the $ \PP $-a.s. sense, and all statements afterwards are to be understood in the $ \widetilde{\PP }$-a.s. sense.\\ \\
\noindent{\bf Step 1. Proof of \eqref{mainthm: mp.spde}.}
\sloppy Since $C^\infty_c(\R)$ is separable under the norm $f \mapsto \max(\Vert f \Vert_\infty, \Vert f' \Vert_\infty, \Vert f''\Vert_\infty)$ (see, e.g., \cite[Lemma 6.1]{lacker2022superposition}), it is enough to show that \eqref{eq: weak sol SPDE} holds for a dense $\{f_i\}_{i \in \N} \subseteq C^\infty_c(\R)$.  To  this end, we follow the strategy in \cite[proof of Proposition 5.4.6]{karatzas2012brownian}. For each $i \in \N$, choose  $ k = 1  $ and $ \phi(x) = x$ on $ \big[-\int_\R |f_i(x) | \, \dd x, \int_\R |f_i(x) |\, \dd  x\big]$ in  \eqref{eq: MRT condition}   to see that 
\begin{align*} 
M^{i}_t :=& \langle u(t,\cdot), f_i\rangle  
- 	 \langle u(0,\cdot), f_i  \rangle  - 	\int_0^t       \Big(\big\langle \Beta(u(r,\cdot)), f_i'\big\rangle + \big\langle (\Sigma + \Gamma)(u(r,\cdot)), f_{i}''\big\rangle \Big) \,  \dd r 
\end{align*}
is an $\FF$-martingale. Similarly, for each $i,j \in \N$,  choose $ k = 2 $ and  $ \phi(x,y)  = xy$ on $ \big[-\int_\R |f_i(x) | \, \dd x, \int_\R |f_i(x) |\, \dd  x\big] \times  \big[-\int_\R |f_j(x) | \, \dd x, \int_\R |f_j(x) |\, \dd  x\big]  $   in  \eqref{eq: MRT condition}  to see that the cross variation  between $M^{i}$ and $M^{j}$ is given by 
\begin{align*}
\left\langle M^{i}, M^{j} \right\rangle_t = \int_0^t  v^i_r \, v^j_r \, \dd r, \quad \text{where } v^i_t := \left\langle G(u(t,\cdot)), f_i' \right\rangle, \,  i \in \N.
\end{align*}
Thus, \eqref{eq: weak sol SPDE} boils down to an extension of the Martingale Representation Theorem, see, e.g., \cite[Proposition 3.4.2]{karatzas2012brownian},  for countably many local martingales. We define for $i,j \in \N$,
\begin{align} \label{def:zij}
z^{i,j}_t := z^{j,i}_t := \frac{\dd}{\dd t} 	\left\langle M^{i}, M^{j} \right\rangle_t = v^i_t \, v^j_t.
\end{align}

For each $ d \in \N $, define the  $d\times d$  matrix-valued process $Z^{(d)}_t := (z^{i,j}_t)_{i, j =1}^d = v^{(d)}_t (v_t^{(d)})^\top $,  where $ v_t^{(d)} := (v^1_t, \dots, v^d_t)^\top $. 
Diagonalizing $Z^{(d)}_t$, we find $d\times d$ matrix-valued processes $ Q^{(d)}_t =  (q_t^{d, i, j})_{i,j = 1}^d $ and $\Lambda^{(d)}_t$
   such that $ (Q^{(d)}_t)^\top Q^{(d)}_t = Id $, 
and $ (Q^{(d)}_t)^\top Z^{(d)}_t Q^{(d)}_t = \Lambda^{(d)}_t $ is diagonal. Moreover, since the rank of $ Z^{(d)}_t $ is at most one, we can assume that $ \Lambda^{(d)}_t $ has  $ |v_t^{(d)}|^2 $ in its $ (1,1) $ entry and  zeros in all other entries, and that the first column of $ Q^{(d)}_t $  is either $v_t^{(d)}/|v_t^{(d)}|$ if $ v^{(d)}_t \neq 0 $ or $ (1,0, \dots, 0)^\top $ otherwise. 

\sloppy  Note that $|q^{d, i,1}_t|\le 1 $ for all $ i=1, \dots, d $, so we can define the $\FF$-martingale
\begin{align}\label{def:Nd}
N^{d}_t = \sum_{i=1}^d\int_0^t q^{d, i,1}_r \, \dd M^{i}_r. 
\end{align}
Now, we can construct an extension $ (\widetilde{\Omega}, \widetilde{\F},\widetilde{\FF} = (\widetilde{\F}_t)_{t \in [0,T]}, \widetilde{\PP}) $ of the probability space $ (\Omega, \F, \FF = (\F_t)_{t \in [0,T]}, \PP) $ supporting independent one-dimensional standard Brownian motions $ \{B^{d}\}_{d \in \N} $ such that each $ B^{d} $ is independent of $ \{N^d\}_{d \in\N} $.  Following the steps of \cite[proof of Proposition 3.4.2]{karatzas2012brownian}, we see that  the process
\begin{align} \label{def:W(d)}
	\widetilde{W}^{d}_t := \int_0^t \mathbf{1}_{\{|v^{(d)}_r| > 0\}}\frac{1}{|v^{(d)}_r|} \, \dd N^{d}_r + \int_0^t  \mathbf{1}_{\{|v^{(d)}_r| = 0\}} 
	\, \dd B^{d}_r
\end{align}
is an $\widetilde{\FF}$-Brownian motion, and $\widetilde{\PP}$-a.s., 
\begin{align*}
M^{i}_t = \int_0^t  q^{d,i,1}_r |v^{(d)}_r| \,\dd \widetilde{W}^{d}_r = \int_0^t v^i_r \, \dd \widetilde{W}^{d}_r, \quad i=1,\dots, d.
\end{align*}

It remains to show that the Brownian motions $\widetilde{W}^{d}$ are the same for all $d \in \N$.  To see this, note that for any $d_1, d_2  \in \N$ and $0  \le t \le T$,  we may deduce  from \eqref{def:W(d)}, \eqref{def:Nd}, \eqref{def:zij} that
\begin{align*}
 \big< \widetilde{W}^{d_1}, \widetilde{W}^{d_2} \big>_t = t.
\end{align*}
Therefore, $\big(\widetilde{W}^{d_1}_t - \widetilde{W}^{d_2}_t\big)^2_{t \in [0,T]}$ is an $ \widetilde{\FF} $-martingale, and thus 
  $\widetilde{\E} \big[(\widetilde{W}^{d_1}_t - \widetilde{W}^{d_2}_t)^2\big]$ = 0 for all $t \in [0,T]$. This shows that 
  $\widetilde{W}^{d_1}$ and $\widetilde{W}^{d_2}$ are indistinguishable, as desired.

\subsection*{Step 2. Proof of $ \E[\W_1(\nu_0, \nu_T)] < \infty$} Next, we show that $ \E[\W_1(\nu_0, \nu_T)] < \infty$. Note that 
\begin{align}
& \E[\W_1(\nu_0, \nu_T)] = \E \left[\int_\R \big|  u(T,x)- u(0,x) \big| \,\dd x \right] \nonumber\\
& \le  \E \left[\int_0^\infty 1 - u(T,x) \, \dd x  \right]  +\E \left[\int_{-\infty}^0 u(T,x) \, \dd x \right] + \int_0^\infty 1 - u(0,x) \,  \dd x    + \int_{-\infty}^0  u(0, x) \,  \dd x.  \label{eq: four terms}
\end{align}
We claim that all four terms in \eqref{eq: four terms} are finite. The third and fourth terms in \eqref{eq: four terms} are finite since $ \nu_0 \in \P_1(\R) $. For the second term,   let $f$ be a smooth, non-increasing function such that $f(x) = 1$ on $(-\infty, 0]$ and $f(x) = 0$ on $[1, \infty)$. Then
\begin{align*}
	  \int_{-\infty}^0  u(T,x) \, \dd x \le \int_\R  u(T,x) \, f(x) \,\dd x. 
\end{align*}
A straightforward approximation argument shows that \eqref{eq: weak sol SPDE} applies to the function $f$ as described. Taking the expectation, we have 
\begin{equation}  \label{eq: expect.u.f.}
	\begin{split}
		& \E \left[\int_{\R}  u(T,x) \, f(x) \, \dd x \right] =    \int_{\R}  u(0,x) \, f(x) \,\dd x   +  \E \left[\int_0^T \! \int_{\R} G(u(r,x)) \, f'(x) \, \dd x \, \dd W_r\right]\\
		&  \hspace{1.0cm}+ \E \left[\int_0^T \!\int_{\R} \Beta(u(r,x)) \, f'(x) +   \Sigma(u(r,x))  \, f''(x) + \Gamma(u(r,x)) \, f''(x) \,  \dd x \, \dd r\right]. 
	\end{split}
\end{equation} 
From Assumption \ref{assum} and the fact that $ f' $ and $ f'' $ are supported on $ [0,1] $, we see that the third term on the right-hand side (RHS) of \eqref{eq: expect.u.f.} is finite. The second term on the RHS of \eqref{eq: expect.u.f.} is zero, as the $\dd W_r$-integrand is bounded.
 Finally, the first term on the RHS of \eqref{eq: expect.u.f.} can be bounded by 
\begin{align*}
	 \int_{\R}  u(0,x) \,\dd x \le \int_{-\infty}^1 u(0, x) \, \dd x,
\end{align*}
which is finite as $ \nu_0 \in \P_1(\R)$. All in all, we deduce that the second term  on the RHS of \eqref{eq: four terms} is finite. The same reasoning applies to the first term as well, completing the proof of $ \E[\W_1(\nu_0, \nu_T)] < \infty$.

In addition, we note that by Young's convolution inequality, for any $\varepsilon > 0$, 
\begin{align} \label{eq: L1.int.u.eps}
	\E \left[\int_\R \big|  u^\varepsilon(T,x)- u^\varepsilon(0,x) \big| \,\dd x \right]  \le \E \left[\int_\R \big|  u(T,x)- u(0,x) \big| \,\dd x \right]< \infty.
\end{align}
\subsection*{Step 3. Proof of \eqref{mainthm: ac} and \eqref{mainthm: u_x in L2}} 
We  apply the definition of weak solution \eqref{eq: weak sol SPDE} to the choice $ f = \varphi_\varepsilon(x - \cdot) $. This is possible by Lemma \ref{lem:moll.com}. As a result, for each $ \varepsilon >0 $ and  $ x  \in\R$, 
\begin{align} \label{eq: u eps dym multi}
	\dd u^\varepsilon (t,x) =\left(-  \Beta(u)^\varepsilon_{x} +  \Sigma(u)^\varepsilon_{xx}  + \Gamma(u)^\varepsilon_{xx} \right) (t,x) \,\dd t - G(u)^\varepsilon_{x}(t,x) \, \dd W_t.
\end{align}
 Note that after mollification, each $ u^\varepsilon(\cdot, x) $ is a semimartingale.
 Using It\^{o}'s Lemma, we deduce 
\begin{align} \label{eq: u.eps.square multi}
\frac{1}{2}\big(u^\varepsilon(T,\cdot)^2  - u^\varepsilon(0,\cdot)^2 \big) =& \int_0^T  u^\varepsilon \big( -  \Beta(u)^\varepsilon_{x}  +  \Sigma(u)^\varepsilon_{xx} + \Gamma(u)^\varepsilon_{xx} \big)  \, \dd t \\
 &+ \frac{1}{2} \int_0^T     \left(G(u)^\varepsilon_{x} \right)^2  \, \dd t - \int_0^T  u^\varepsilon  \, G(u)^\varepsilon_x \, \dd W_t. \label{eq: stoc int term}
\end{align}
	We claim that the stochastic integral  is a martingale. Note that the measure $\dd G(u(t, \cdot))$ is finite, as \cite[Theorem 31.2]{billingsley2017probability} implies
\begin{align*}
 \int_\R \dd G(u(t, \cdot)) \le G(1) - G(0) = \int_0^1 \gamma(a)  \,\dd a< \infty. 
\end{align*} 
Thus, by Jensen's inequality,
\begin{align*}
	\int_0^T \big(G(u)^\varepsilon_x\big)^2 \, \dd t \le \big(G(1) - G(0)\big) \int_0^T \!\int_\R \varphi_\varepsilon^2(\cdot- y) \, \dd G(u(t,y)) \, \dd t \le \big(G(1) - G(0)\big)^2 \, \Vert \varphi_\varepsilon   \Vert^2_\infty \, T.
\end{align*}
Therefore, \cite[Corollary IV.1.25]{revuz2013continuous} implies that the stochastic integral in \eqref{eq: stoc int term} is a  martingale. Taking the expectation  in \eqref{eq: u.eps.square multi}--\eqref{eq: stoc int term}, we have
\begin{align*}
\frac{1}{2} \E \left[u^\varepsilon(T,\cdot)^2 - u^\varepsilon(0,\cdot)^2\right]   = & \E \left[\int_0^T  u^\varepsilon \big( -  \Beta(u)^\varepsilon_{x}  +  \Sigma (u)^\varepsilon_{xx} + \Gamma(u)^\varepsilon_{xx} \big) + \frac{1}{2} \left(G(u)^\varepsilon_{x} \right)^2 \dd t \,  \right].
\end{align*}

Note that  since each $ u(t, \cdot) $ is a CDF,  so is $ u^\varepsilon(t, \cdot) $. In conjunction with \eqref{eq: L1.int.u.eps},
\begin{align} \label{DCT.bound}
	\E\left[\int_\R \big| u^\varepsilon(T,x)^2 - u^\varepsilon(0, x)^2 \big| \, \dd x \right]\le 2 \E \left[\int_\R \big| u^\varepsilon(T,x) - u^\varepsilon(0,x) \big| \, \dd x \right] < \infty. 
\end{align}
 Also, the convolution of a finite measure with the Gaussian kernel or its derivatives is in $ L^p $ for any $ p \in [1,\infty) $. 
Hence, integrating over $ x $ in $ \R $ and using Fubini's Theorem,  we have
\begin{align} \label{eq: E.u.eps.square}
\begin{split}
&\frac{1}{2}   \E \left[ \int_\R u^\varepsilon(T,x)^2 - u^\varepsilon(0,x)^2 \, \dd x\right]   \\
& =   \E \left[ \int_0^T \!\int_\R u^\varepsilon \big( -  \Beta(u)^\varepsilon_{x}  +  \Sigma(u)^\varepsilon_{xx} + \Gamma(u)^\varepsilon_{xx} \big) \, \dd x \, \dd t + \frac{1}{2}  \int_0^T \! \int_\R \left(G(u)^\varepsilon_{x} \right)^2 \, \dd x \, \dd t\right].
\end{split}
\end{align}

\noindent{\bf Step 3.1. Convergence of LHS.}
We take $\varepsilon \downarrow 0$ on the left-hand side (LHS) of \eqref{eq: E.u.eps.square}.  Note that
\begin{align}
&\left|  \E \left[\int_{\R} u^\varepsilon(T,x)^2 - u^\varepsilon(0,x)^2 \, \dd x\right]      -  \E \left[ \int_{\R} u(T,x)^2 - u(0,x)^2 \, \dd x  \right]  \right| \nonumber \\
& \le   \E \left[\int_{\R} \bigg| \Big(\big(u^\varepsilon(T,x) - u^\varepsilon(0,x)\big) - \big(u(T,x) - u(0,x)\big)\Big) \big(u^\varepsilon(T,x) + u^\varepsilon(0,x)\big) \bigg| \, \dd x   \right] \label{eq: first integral multi}\\
&  \hspace{0.4cm}+ \E \left[\int_{\R} \bigg| \Big(\big(u^\varepsilon(T,x) + u^\varepsilon(0,x)\big)- \big(u(T,x)+ u(0,x)\big)\Big)  \big(u(T,x)- u(0,x)\big) \bigg| \,\dd x\right].   \label{eq: second integral multi}
\end{align}
	Let us study the term in \eqref{eq: second integral multi} first. As $ \varepsilon \downarrow 0 $, the integrand converges to $ 0 $ for Lebesgue-a.e. $x \in \R$ by Lemma \ref{lem: CDF}(i).   Also, the  integrand is dominated by $ 2 \, |u(T, \cdot) - u(0, \cdot)|$, which is in $ L^1(\Omega \times \R) $ by Step 2. Thus the Dominated Convergence Theorem implies that \eqref{eq: second integral multi}  converges to $0$ as $\varepsilon \downarrow 0$. 
	
	Turning to \eqref{eq: first integral multi}, it is bounded by
	\begin{align} \label{eq: first.integral.bound}
	2 \, \E \left[\int_{\R} \Big| \big(u^\varepsilon(T,x) - u^\varepsilon(0,x)\big) - \big(u(T,x) - u(0,x)\big) \Big| \, \dd x \right].
	\end{align} 
	We know that  $ u(T,\cdot) - u(0,\cdot) \in L^1(\R)$ a.s. by Step 2. Therefore, Lemma \ref{lem: CDF}(ii) implies  that $ u^\varepsilon(T,\cdot) - u^\varepsilon(0,\cdot) $ converges to $ u(T,\cdot) - u(0,\cdot)$ in $ L^1(\R) $ a.s.
Moreover, by Young's convolution inequality, 
\begin{align*}
	\big\Vert \big(u(T, \cdot) - u(0, \cdot)\big)^\varepsilon \big\Vert_{L^1(\R)}  \le \Vert \varphi_\varepsilon \Vert_{L^1(\R)} \Vert u(T, \cdot) - u(0, \cdot)  \Vert_{L^1(\R)} = \W_1(\nu_0, \nu_T), 
\end{align*}
and so the term inside the expectation of \eqref{eq: first.integral.bound} is bounded by $ 2 \, \W_1(\nu_0, \nu_T) $. Together with  $ \E[\W_1(\nu_0, \nu_T)] < \infty$ from Step 2, the Dominated Convergence Theorem implies that \eqref{eq: first.integral.bound} converges to $ 0 $ as $ \varepsilon \downarrow 0 $. Hence, 
\begin{align*}
	\lim_{\varepsilon \downarrow 0} \E \left[ \int_{\R}  u^\varepsilon(T,x)^2 - u^\varepsilon(0,x)^2 \,  \dd x  \right]  =  \E \left[ \int_{\R} u(T,x)^2 - u(0,x)^2 \,  \dd x\right].  
\end{align*}

			\noindent{\bf Step 3.2. $\left( u^\varepsilon_{x}\right)_{\varepsilon \downarrow 0}$ is weakly compact in $L^2([0,T] \times \R)$ a.s.} 
	From Step 3.1, we know 
		\begin{align*} 
		&\lim_{\varepsilon \downarrow 0}    \E \left[ \int_0^T \!  \int_{\R} u^\varepsilon \big(-  \Beta(u)^\varepsilon_{x}  +  \Sigma(u)^\varepsilon_{xx}  + \Gamma(u)^\varepsilon_{xx}  \big)  \, \dd x \,  \dd t  + \frac{1}{2}  \left(G(u)^\varepsilon_{x} \right)^2\, \dd x\, \dd t \right]  \\
			&= 	\frac{1}{2} \, \E \left[ \int_{\R} u(T,x)^2 - u(0,x)^2 \,  \dd x\right]. 
		\end{align*} 
		In a manner similar to \eqref{DCT.bound}, we see that the RHS is finite.
		Thus
		\begin{align} \label{eq:finitness2.25}
\lim_{\varepsilon \downarrow 0}   \E \left[ \int_0^T \! \int_{\R} u^\varepsilon \big(-  \Beta(u)^\varepsilon_{x}  +  \Sigma(u)^\varepsilon_{xx}  + \Gamma(u)^\varepsilon_{xx}  \big)   + \frac{1}{2} \left(G(u)^\varepsilon_{x} \right)^2 \, \dd x\, \dd t \right]     \in \R.
		\end{align}  
		
		We now show that 
		\begin{align} \label{eq:finiteness2.26}
			\limsup_{\varepsilon \downarrow 0}   \E \left[ \int_0^T \! \int_{\R}  u^\varepsilon  \,  \Gamma(u)^\varepsilon_{xx}   + \frac{1}{2}  \left(G(u)^\varepsilon_{x} \right)^2   \dd x \,  \dd t    \right] < \infty .
		\end{align}
		To see this, first fix $ \varepsilon > 0 $,  $t \in [0,T]$ and $ -\infty < a^- < a^+ < \infty $. An integration by parts gives
		\begin{align*}
			\int_{a^-}^{a^+}  u^\varepsilon(t,x) \,  \Gamma(u)^\varepsilon_{xx}(t, x)\,   \dd x = u^\varepsilon(t, \cdot)  \, \Gamma(u)^\varepsilon_x(t,\cdot) \Big|_{a^-}^{a^+}  - \int_{a^-}^{a^+}   u^\varepsilon_x(t,x) \, \Gamma(u)^\varepsilon_x (t,x) \,  \dd x. 
		\end{align*}
		We claim that the boundary terms vanish as $ a^- \to -\infty $ and $ a^+ \to \infty $ along suitable sequences $a^-_k \to -\infty$ and $a^+_k \to \infty$, respectively. If this were not the case, this would imply $ u^\varepsilon(t, x)  \, \Gamma(u)^\varepsilon_x(t,x) $ is bounded away from zero for $ |x| $ sufficiently large.  Then, $ \int_\R u^\varepsilon(t, x)  \, \Gamma(u)^\varepsilon_x(t,x)  \,  \dd x = \infty $, contradicting the fact that 
		 \begin{align*}
		  \int_\R u^\varepsilon(t, x)  \, \Gamma(u)^\varepsilon_x(t,x)  \,  \dd x \le \int_\R  \Gamma(u)^\varepsilon_x(t,x)  \,  \dd x < \infty
		 \end{align*}
		 due to the finiteness of the measure $  \Gamma(u)^\varepsilon_x(t,x)  \,  \dd x  $. 
	Therefore, the quantity in \eqref{eq:finiteness2.26} equals to 
	\begin{align} \label{nonpos}
		\limsup_{\varepsilon \downarrow 0} \E \left[ \int_0^T \!\int_\R -  u^\varepsilon_x \,  \Gamma(u)^\varepsilon_x + \frac{1}{2}  \big(G(u)^\varepsilon_x\big)^2 \, \dd x \,\dd t \right].
	\end{align}
	Let $u^{-1}(t,\xi)$ be the $\xi$-quantile of $\dd u (t, \cdot)$. By Fubini's theorem, the integrand is
	\begin{align*}
		&-  u^\varepsilon_x \,  \Gamma(u)^\varepsilon_x (t,x) + \frac{1}{2}  \big(G(u)^\varepsilon_x\big)^2 (t,x) \\
		&= - \frac{1}{2} \int_\R \varphi_\varepsilon'(y) \int_0^{u(t, x-y)} 1\, \dd \xi \, \dd y \cdot \int_\R \varphi_{\varepsilon}'(y) \int_0^{u(t,x-y)} \gamma^2(\xi)\, \dd \xi \, \dd y \\
		& \hspace{0.45cm} + \frac{1}{2} \bigg(\int_\R \varphi_\varepsilon'(y) \int_0^{u(t,x-y)} \gamma(\xi) \, \dd \xi \, \dd y\bigg)^2 \\
		&= - \frac{1}{2}\int_0^1 \varphi_{\varepsilon}\big(x - u^{-1}(t, \xi)\big) \, \dd \xi \cdot  \int_0^1 \varphi_\varepsilon \big(x - u^{-1}(t, \xi) \big) \gamma^2(\xi) \, \dd \xi \\
		&\hspace{0.45cm} + \frac{1}{2} \bigg(\int_\R \varphi_\varepsilon\big(x - u^{-1}(t, \xi)\big) \, \gamma(\xi) \, \dd \xi \bigg)^2
	\end{align*}
	which is non-positive by the Cauchy-Schwarz inequality. Thus \eqref{nonpos} is also non-positive.

		Together with \eqref{eq:finitness2.25}, we deduce that
		\begin{align*}
			\liminf_{\varepsilon \downarrow 0} \E \left[\int_0^T \! \int_{\R} u^\varepsilon \big(- \Beta(u)^\varepsilon_{x}   + \Sigma(u)^\varepsilon_{xx} \big) \,  \dd x \, \dd t \right]  > -\infty.
		\end{align*}
		Note further that for fixed $ \varepsilon > 0 $,  $t \in [0,T]$ and $ -\infty < a^- < a^+ < \infty $, an integration by parts gives
		\begin{align*}
			-\int_{a^-}^{a^+} u^\varepsilon(t,x) \,  \Beta(u)^\varepsilon_{x}(t,x) \,  \dd x & = -  u^\varepsilon(t, \cdot) \, \Beta(u(t, \cdot))^\varepsilon \Big|_{a^-}^{a^+} +   \int_{a^-}^{a^+} u^\varepsilon_{x}(t,x)\,  \Beta(u)^\varepsilon (t,x)\,  \dd x  \\
			&\le 2 \, \Vert \Beta \Vert_\infty  + \Vert \Beta \Vert_\infty \int_{a^-}^{a^+} u^\varepsilon_{x}(t,x)  \,  \dd x  \le 3 \, \Vert \Beta \Vert_\infty,
		\end{align*}
		and therefore,
		\begin{align*}
			- \int_0^T \int_\R u^\varepsilon \,  \Beta(u)^\varepsilon_{x} \,  \dd x \, \dd t  \le 3 \, \Vert \Beta \Vert_\infty  \, T. 
		\end{align*}
		This implies 
			\begin{align*}
			\liminf_{\varepsilon \downarrow 0}   \E \left[\int_0^T \int_{\R} u^\varepsilon  \,  \Sigma(u)^\varepsilon_{xx} \,  \dd x \,  \dd t\right] > -\infty.
		\end{align*}
		
		An integration-by-parts argument as above shows that 
		\begin{align*} 
			\int_\R u^\varepsilon \,  \Sigma(u)^\varepsilon_{xx} \,  \dd x =   - \int_\R u^\varepsilon_x \, \Sigma(u)^\varepsilon_x \,  \dd x ,
		\end{align*}
		and so we deduce  that
		\begin{align} \label{eq: finiteness.2.21}
			\limsup_{\varepsilon \downarrow 0}    \E \left[ \int_0^T \! \int_{\R} u^\varepsilon_{x} \,  \Sigma(u)^\varepsilon_x \,  \dd x \,  \dd t\right] < \infty.
		\end{align}
		In view of Fubini's Theorem and Assumption \ref{assum}\eqref{eq: assum sigma},  we have
		\begin{align*}
			\Sigma(u)^\varepsilon_x (t,x) &=   \frac{1}{2} \int_\R \int_0^{u(t, x - y)} \sigma^2(\xi) \, \dd \xi \, \varphi'(y) \, \dd y  \\
			&= \frac{1}{2} \int_0^1 \varphi_\varepsilon \big(x - u^{-1}(t, \xi) \big) \, \sigma^2(\xi) \, \dd \xi  \ge \frac{c_\sigma^2}{2} \, u^\varepsilon_x (t,x),
		\end{align*}
		where $c _\sigma := \inf_{a \in [0,1]} \sigma(a) > 0$. Together with \eqref{eq: finiteness.2.21}, we see that
		\begin{align*}
			&\limsup_{\varepsilon \downarrow 0}   \E \left[\int_0^T\! \int_{\R} \big(u^\varepsilon_{x}\big)^2 \,  \dd x \,  \dd t \right]
			< \infty.
		\end{align*}
		By Fatou's Lemma, 
		\begin{align*}
			\E \left[\liminf_{\varepsilon \downarrow 0} \int_0^T\! \int_{\R} \big(u^\varepsilon_{x}\big)^2 \,  \dd x \,  \dd t \right] \le  \liminf_{\varepsilon \downarrow 0} \E \left[\int_0^T \! \int_{\R} \big(u^\varepsilon_{x}\big)^2 \,  \dd x \,  \dd t  \right] < \infty,
		\end{align*}
		which implies $ \PP $-a.s.,
		\begin{align*}
			\liminf_{\varepsilon \downarrow 0} \int_0^T \! \int_{\R} \big(u^\varepsilon_{x}\big)^2 \,  \dd x \,  \dd t  < \infty. 
		\end{align*}
		By the Banach-Alaoglu Theorem, there exists a (random) subsequence $\left( u^{\varepsilon_n}_{x}\right)_{n \in \N}$ and a unique  $v \in L^2([0,T] \times \R)$ such that for all $g \in L^2([0,T] \times \R)$,
		\begin{align} \label{eq: conv 1 multi}
			\lim_{n \to \infty} \int_0^T \! \int_{\R} g \,  u^{\varepsilon_n}_{x} \,  \dd x \,  \dd t  = \int_0^T\!\int_{\R} g \, v \,  \dd x \,  \dd t.
		\end{align}
			\noindent{\bf Step 3.3. Completing the proof.}
		To conclude Step 3, we note that for any $g \in C^\infty_c ([0,T] \times \R)$ and $ t \in [0,T] $, $\int_\R g(t,x) \, u^\varepsilon_x(t,x)\, \dd x \stackrel{\varepsilon \downarrow 0 }{\to} \int_\R g(t,x) \, \nu_t(\dd x) $. Thus, 
		the Bounded Convergence Theorem implies that for any $g \in C^\infty_c ([0,T] \times \R)$,
		\begin{align}\label{eq: conv 2 multi}
			\lim_{\varepsilon \downarrow 0} \int_0^T \! \int_{\R}  g(t, x) \,  u^\varepsilon_{x}(t, x) \,  \dd x \,  \dd t = 	\int_0^T \! \int_{\R}  g(t, x) \,  \nu_t(\dd x) \,  \dd t.  
		\end{align}
		Since $C^\infty_c([0,T] \times \R)$ is a distribution-determining class, comparing \eqref{eq: conv 1 multi} and \eqref{eq: conv 2 multi} shows that  $   \nu_t(\dd x)  \, \dd t  $ has  density $v \in L^2([0,T] \times \R)$, which yields part \eqref{mainthm: u_x in L2} of the theorem.  In addition, for Lebesgue a.e. $t \in [0, T]$, $ \nu_t(\dd x) $ is absolutely continuous, which shows part \eqref{mainthm: ac} of the theorem.
\subsection*{Step 4}  The next two steps are preparations for the proof of part  \eqref{mainthm: pe} of the theorem. Fix $ 0 \le s \le t \le T $ and $ (\eta, y) \in \R^2$. To lighten notation, we  abbreviate $ \varrho (\xi, t, x; \eta ,y)$ defined in \eqref{eq: varrho multi} by $ \varrho (\xi, t, x)$.  In this step, we  show the key identity
\begin{align} \label{eq: Step 1 result}
	\sum_{i=1}^{11} I_i= I_{12} + I_{13},
\end{align}
where with $ u:= u(r,x) $:
\begin{align*}
	\begin{split}
			I_1 &:= \int_s^t  \int_{\R}  \int_0^u \varrho^0_{x}(x - y - b(\xi) r  - \gamma(\xi)  W_r, \xi - \eta) \, b(\xi)\, \dd \xi \,  \dd x \, \dd r, \\
		I_2&: = \int_s^t \int_{\R} \varrho^0(x  - y - b(u)r - \gamma(u)  W_r, u - \eta) \,
		\Beta(u)_{x} \, \dd x \, \dd r, \\
		I_3 &:= - \int_s^t  \int_{\R}  \varrho^0(x - y - b(u) r  - \gamma(u)  W_r, u - \eta)  \, \Sigma(u)_{xx}  \, \, \dd x \, \dd r, \\
		I_4 &:= - \int_s^t  \int_{\R}  \varrho^0(x- y - b(u) r  - \gamma(u)  W_r , u -\eta)  \, \Gamma(u)_{xx}  \,  \dd x \, \dd r, \\
		I_5 &:=  \int_s^t \int_{\R}  \varrho^0(x -y - b(u) r  - \gamma(u)  W_r, u - \eta) \, G(u)_{x} \,\dd x \, \dd W_r, \\
		I_6 &:= \frac{1}{2}  \int_s^t \int_{\R}  \varrho^0_{x}(x - y- b(u) r  - \gamma(u) W_r , u - \eta) \, \big(b'(u)r + \gamma'(u) \, W_r \big) \left(G(u)_{x}\right)^2 \dd x \, \dd r, \\
		I_7 &:= - \frac{1}{2} \int_s^t \int_{\R}  \varrho^0_\xi(x - y - b(u) r  - \gamma(u) W_r , u - \eta)  \left(G(u)_{x}\right)^2  \dd x \, \dd r, \\
		I_8 &:= \int_s^t \int_{\R} \int_0^u \varrho^0_{x}(x - y - b(\xi) r - \gamma(\xi) W_r , \xi - \eta) \, \gamma(\xi) \, \dd \xi \, \dd x  \, \dd W_r,\\
		I_9 &:=  -\frac{1}{2} \int_s^t \int_{\R}\int_0^u \varrho^0_{xx}(x  - y - b(\xi)r - \gamma(\xi) W_r, \xi - \eta) \, \gamma^2(\xi) \, \dd \xi \,\dd x  \,\dd r, \\
		I_{10} &:= - \int_s^t \int_{\R}\varrho^0_{x}(x - y - b(u) r - \gamma(u)  W_r , u - \eta) \, \gamma(u)\, G(u)_x \,\dd x \, \dd r,\\
		I_{11} &:=  \frac{1}{2}\int_s^t \int_{\R} \int_0^u  \varrho^0_{xx}(x - y - b(\xi) r- \gamma(\xi) W_r , \xi -\eta) \, \sigma^2(\xi) \,\dd \xi \,\dd x \, \dd r,\\
		I_{12} &:= - \frac{1}{2} \int_s^t \int_{\R}  \varrho_{x}^0 \left(x - y - b(u) r- \gamma(u)  W_r , u - \eta\right) \left(b'(u) r + \gamma'(u) W_r \right)  \sigma^2(u) \, u^2_{x} \, \dd x \,\dd r,\\
		I_{13} &:=  \frac{1}{2} \int_s^t \int_{\R}  \varrho_\xi^0 (x - y - b(u) r - \gamma(u)  W_r , u - \eta)\,\sigma^2(u) \, u^2_{x} \, \dd x \, \dd r.   
	\end{split}
\end{align*}
\noindent{\bf Step 4.1.}
We first show that $I_1+ I_2 = 0$ and $I_5 + I_8 = 0$. Since they are similar, we  focus on  $I_5 + I_8 = 0$. We make a change of variables in $ I_5 $. Note that $ G(u)_x = \gamma(u) \, u_x $ for a.e. $r \in [s,t]$. We claim that for a.e. $r \in [s,t]$, 
\begin{equation}  
	\label{eq: co-area claim 1 multi}
	\begin{split}
			&\int_{\R}  \varrho^0\left(x - y - b(u) r - \gamma(u) W_r , u- \eta \right) \gamma(u) \,u_{x} \,\dd x \\&=  \int_0^1 \varrho^0 \left(  u^{-1}(r,  \xi) - y - b(\xi) r - \gamma(\xi)W_r , \xi - \eta\right) \gamma(\xi) \, \dd \xi,
	\end{split}
\end{equation}
where $u^{-1}(r, \xi)$ is the $\xi$-quantile of $\dd u(r, \cdot)$.  Indeed, from part \eqref{mainthm: ac} of the theorem, we know that $ u(r,\cdot) $ is continuous for  a.e. $r \in [s,t]$. For any such $ r \in [s,t] $, the identity \eqref{eq: co-area claim 1 multi}  follows from the co-area formula of Fleming-Rishel  \cite[Theorem 1]{fleming1960integral} (see also  \cite[equation (1.4)]{dierkes2023remarks}, where we take $f(x) := u(r,x) $ and  $ g(x) := \varrho^0\left(x -y - b(u) r -  \gamma(u) W_r, u - \eta\right) \, \gamma(u)   $).  Hence,
\begin{align*}
	I_5 &=  \int_s^t \int_0^1 \varrho^0 \big(u^{-1}(r,\xi) - y -b(\xi) r -  \gamma(\xi) W_r ,  \xi - \eta\big) \gamma(\xi) \, \dd \xi \, \dd W_r.  
\end{align*}

On the other hand, by Fubini's theorem,
\begin{align}
	I_8 &= \int_s^t \int_0^1  \int_\R  \mathbf{1}_{\{0 \le \xi \le u(r, x)\}} \, \varrho^0_{x}(x - y - b(\xi)r - \gamma(\xi) W_r , \xi - \eta) \, \dd x \, \gamma(\xi) \, \dd \xi  \, \dd W_r. \label{eq: (f) intermediate}
\end{align}
 For any fixed  $ r \in [s,t] $ and for all $ \xi \in (0,1)$, 
the innermost integral is
\begin{align*}
	 \int_{u^{-1}(r, \xi)}^\infty \varrho^0_{x}(x- y - b(\xi)r - \gamma(\xi)  W_r , \xi - \eta) \, \dd x 
= - \varrho^0\left( u^{-1}(r, \xi) - y - b(\xi) r - \gamma(\xi)W_r ,\xi - \eta\right).
\end{align*}
Putting this back into \eqref{eq: (f) intermediate}, we see that $I_5 + I_8 = 0$ as claimed.

\noindent{\bf Step 4.2.}
From Step 4.1, we know $I_1 + I_2 = 0$ and $I_5 + I_8 = 0$. So, we are left to show 
\begin{align} \label{eq: key.identity.intermediate}
I_3 + I_4 + I_6 + I_7 + I_9 + I_{10} + I_{11} = I_{12} + I_{13}.
\end{align}
 Integrating by parts, and noting that $\Sigma(u)_{x} = \frac{1}{2} \, \sigma^2 (u) \, u_{x}$, we have
\begin{equation}  \label{eq: term 1(b) integrated multi}
	\begin{split}
			I_3 =   \frac{1}{2}\int_s^t \int_{\R}\Big( & \varrho^0_{x}\left(x- y - b(u)r -  \gamma(u) W_r , u -\eta\right) 
\\
&-  \varrho^0_{x}\left(x - y - b(u)r  - \gamma(u)  W_r , u - \eta\right)\left(b'(u) r + \gamma'(u) W_r \right) u_{x} \\
		& + \varrho^0_\xi (x  - y - b(u)r - \gamma(u)  W_r, u - \eta)u_{x}\Big) \sigma^2(u) \, u_{x} \, \dd x \, \dd r \\
		&= I_{14} + I_{12} + I_{13},
	\end{split}
\end{equation}
where 
\begin{align*}
	I_{14} &:=  \frac{1}{2}\int_s^t \int_{\R}  \varrho^0_{x}\left(x - y - b(u)r -  \gamma(u) W_r , u - \eta\right) \,  \sigma^2(u) \,  u_{x} \, \dd x \,  \dd r.
\end{align*}
We claim that $ I_{11} = -I_{14} $. Indeed, using the same change-of-variable technique as in Step 4.1, we see that
\begin{align*}
	I_{14} &= \frac{1}{2}\int_s^t \int_0^1 \varrho^0_{x} \big( u^{-1}(r, \xi)  - y - b(\xi )r - \gamma(\xi)W_r, \xi - \eta\big) \sigma^2(\xi) \,  \dd \xi \,  \dd r.
\end{align*}
Similarly to \eqref{eq: (f) intermediate}, we obtain
\begin{align*}
	I_{11} &=  \frac{1}{2}\int_s^t \int_0^1  \int_{u^{-1}(r, \xi)}^\infty \varrho^0_{xx}(x - y - b(\xi)r - \gamma(\xi)   W_r , \xi - \eta) \, \dd x \, \sigma^2(\xi) \, \dd \xi \,  \dd r =  -I_{14}.
\end{align*}
\noindent{\bf Step 4.3.}
From Step 4.2, we have $ I_3 = I_{14} + I_{12} + I_{13} $ and $ I_{11} = -I_{14} $. Thus, on account of \eqref{eq: key.identity.intermediate},  we are left to show 
\begin{align*} 
	  I_{4} + I_{6} + I_7 + I_9 + I_{10} =  0.
\end{align*}
In the same way as proving  $ I_3 = I_{14} + I_{12} + I_{13} $ in Step 4.2, we can show that $ I_4 =  -\frac{1}{2}I_{10} - I_{6} - I_{7}$. Also, similar to how $ I_{11} = -I_{14} $ was proven in Step 4.2, we can show that $ I_{9} = -\frac{1}{2} I_{10} $. This completes the proof of \eqref{eq: Step 1 result}.
\\
\subsection*{Step 5} 
In this step, we  define mollified versions of the terms $I_1, \dots, I_{10}$ and show their convergences to $I_1, \dots, I_{10}$ as the mollification parameter vanishes. Fix $ 0 \le s \le t \le T $ and $ (\eta, y) \in \R^2$. As in Step 4, we abbreviate $ \varrho (\xi, t, x; \eta ,y)$  in \eqref{eq: varrho multi} by $ \varrho (\xi, t, x)$.  We begin by defining 
\begin{align} \label{a.to.j.epsilon}
	\begin{split} 
		I^\varepsilon_1 &= \int_s^t  \int_{\R}  \int_0^{u^\varepsilon} \varrho^0_{x}(x  - y - b(\xi) r  - \gamma(\xi)  W_r, \xi - \eta) \, b(\xi)\, \dd \xi \,  \dd x \, \dd r, \\
		I^\varepsilon_2& =  \int_s^t \int_{\R} \varrho^0(x -y - b(u^\varepsilon)r - \gamma(u^\varepsilon)  W_r , u^\varepsilon -\eta) \,
		\Beta(u)^\varepsilon_{x} \, \dd x \, \dd r, \\ 
			I^\varepsilon_3 &= -  \int_s^t  \int_{\R}  \varrho^0(x - y - b(u^\varepsilon)r - \gamma(u^\varepsilon)  W_r , u^\varepsilon - \eta)  \, \Sigma(u)^\varepsilon_{xx}  \,  \dd x \, \dd r, \\
		I^\varepsilon_4 &= - \int_s^t  \int_{\R}  \varrho^0(x - y - b(u^\varepsilon)r  - \gamma(u^\varepsilon)  W_r , u^\varepsilon - \eta)  \, \Gamma(u)_{xx}^\varepsilon  \,  \dd x \, \dd r, \\
		I^\varepsilon_5 &= \int_s^t \int_{\R}  \varrho^0(x - y - b(u^\varepsilon)r -  \gamma(u^\varepsilon)  W_r , u^\varepsilon - \eta) \, G(u)_{x}^\varepsilon \, \dd x \, \dd W_r, \\
		I^\varepsilon_6 &= \frac{1}{2} \int_s^t \int_{\R}  \varrho^0_{x}(x - y - b(u^\varepsilon)r  - \gamma(u^\varepsilon)  W_r , u^\varepsilon - \eta) \, \big(b'(u^\varepsilon)r + \gamma'(u^\varepsilon) \, W_r \big) \left(G(u)_{x}^\varepsilon\right)^2\dd x \, \dd r, \\
		I^\varepsilon_7 &= - \frac{1}{2}  \int_s^t \int_{\R}  \varrho^0_\xi(x-y  - b(u^\varepsilon)r  -  \gamma(u^\varepsilon) W_r , u^\varepsilon- \eta)  \left(G(u)_{x}^\varepsilon\right)^2 \dd x \, \dd r, \\
		I^\varepsilon_8 &=  \int_s^t \int_{\R} \int_0^{u^\varepsilon} \varrho^0_{x}(x - y - b(\xi) r - \gamma(\xi) W_r , \xi - \eta) \, \gamma(\xi) \, \dd \xi \, \dd x  \, \dd W_r,\\
		I^\varepsilon_9 &=  -\frac{1}{2} \int_s^t \int_{\R}\int_0^{u^\varepsilon} \varrho^0_{xx}(x - y - b(\xi) r - \gamma(\xi)W_r , \xi - \eta) \, \gamma^2(\xi) \, \dd \xi \, \dd x  \, \dd r, \\
		I^\varepsilon_{10} &= -   \int_s^t \int_{\R}\varrho^0_{x}(x - y - b(u^\varepsilon)r -  \gamma(u^\varepsilon)  W_r, u^\varepsilon -\eta) \, \gamma(u^\varepsilon) \, G(u)_{x}^\varepsilon \,\dd x \, \dd r.
	\end{split}
\end{align}

We aim to show the  convergences 
\begin{align} \label{eq: convergence 1}
	\lim_{\varepsilon \downarrow 0}\int_{\R^{2}} \overline{\chi}(\xi, u^\varepsilon(\cdot ,x)) \,  \varrho(\xi, \cdot, x ) \, \dd \xi \, \dd x \bigg|_s^t = \int_{\R^{2}} \overline{\chi}(\xi, u(\cdot ,x))  \, \varrho(\xi, \cdot,x ) \, \dd \xi \, \dd x \bigg|_s^t 
\end{align}
and 
\begin{align} \label{eq: convergence 2}
	\lim_{\varepsilon \downarrow 0} I^\varepsilon_i = I_i, \quad i =1, \dots, 10. 
\end{align}
We  divide them into three groups. \\ \\
\noindent{\it Group 1: } \eqref{eq: convergence 1}, $ \lim_{\varepsilon \downarrow 0} I^\varepsilon_1 = I_1$,  and $ \lim_{\varepsilon \downarrow 0} I^\varepsilon_9 = I_9$. \\ \\
Let us  prove \eqref{eq: convergence 1} first. Recalling the definitions of $ \overline{\chi } $ in \eqref{chi} and $ \varrho $ in \eqref{eq: varrho multi}, and that $ 0 \le u^\varepsilon \le 1$ because each $ u(t, \cdot) $ is assumed to be  a CDF,  we are left to show that for any fixed $ t \in [0,T] $,
\begin{align*}
	&\lim_{\varepsilon \downarrow 0} \int_\R \int_0^{u^\varepsilon(t,x)} \varrho^0 (x - y- b(\xi)t - \gamma(\xi)W_t , \xi -\eta) \, \dd \xi \, \dd x \\
	&= \int_\R \int_0^{u(t,x)} \varrho^0 (x - y - b(\xi)t - \gamma(\xi)W_t , \xi -\eta) \, \dd \xi \, \dd x. 
\end{align*}
The convergence of the inner integral is a consequence of Lemma \ref{lem: CDF}(i), and the convergence of the outer integral then follows from the Dominated Convergence Theorem thanks to the compact support of $ \varrho^0 $ and the boundedness of $ b $ and $ \gamma $. 

The  other two convergences $ \lim_{\varepsilon \downarrow 0} I^\varepsilon_1 = I_1$  and $ \lim_{\varepsilon \downarrow 0} I^\varepsilon_9 = I_9$ are obtained from the Dominated Convergence Theorem in a similar manner.  \\

\noindent{\it Group 2: }	$ \lim_{\varepsilon \downarrow 0} I^\varepsilon_2 = I_2 $,  $ \lim_{\varepsilon \downarrow 0} I^\varepsilon_3 = I_3 $,  $ \lim_{\varepsilon \downarrow 0} I^\varepsilon_4 = I_4 $, $ \lim_{\varepsilon \downarrow 0} I^\varepsilon_6 = I_6 $, $ \lim_{\varepsilon \downarrow 0} I^\varepsilon_7 = I_7 $,  and $ \lim_{\varepsilon \downarrow 0} I^\varepsilon_{10} = I_{10} $. \\ \\
Since all statements in this group can be proven similarly, we just prove $ \lim_{\varepsilon \downarrow 0} I^\varepsilon_3 = I_3 $. Recall that from  \eqref{eq: term 1(b) integrated multi},  we have the identity $ I_3 = I_{14} + I_{12} + I_{13} $.  Similarly, we can show that $ I^\varepsilon_3 = I^\varepsilon_{14} + I^\varepsilon_{12} + I^\varepsilon_{13} $, where 
\begin{align*}
	I^\varepsilon_{14}  &:=  \frac{1}{2} \int_s^t \int_{\R} \varrho^0_{x}\left(x  - y - b(u^\varepsilon)r -  \gamma(u^\varepsilon)  W_r, u^\varepsilon - \eta\right) \big(\sigma^2(u)u_{x}\big)^\varepsilon \, \dd x \, \dd r, \\
I^\varepsilon_{12} &:= -   \frac{1}{2}\int_s^t \int_{\R} \varrho^0_{x}(x - y- b(u^\varepsilon)r -  \gamma(u^\varepsilon)  W_r , u^\varepsilon -\eta) \left(b'(u^\varepsilon)r +  \gamma'(u^\varepsilon) W_r \right)u^\varepsilon_{x} \big(\sigma^2(u) \, u_{x}\big)^\varepsilon\, \dd x \,\dd r, \\
	I^\varepsilon_{13}  &:=   \frac{1}{2}  \int_s^t \int_{\R}   \varrho^0_\xi (x - y - b(u^\varepsilon)r -  \gamma(u^\varepsilon)  W_r, u^\varepsilon -\eta) \, u^\varepsilon_{x} \big(\sigma^2(u) \, u_{x}\big)^\varepsilon \, \dd x \, \dd r .
\end{align*}

Firstly, we show that $ I^\varepsilon_{14} \to I_{14} $. From  Step 3 and Assumption \ref{assum}\eqref{eq: assum sigma}, we know that  $ \sigma^2(u) \, u_{x} \in L^2([0,T] \times \R) $.  We claim that 
\begin{align}
	(\sigma^2(u)\,u_{x})^\varepsilon \to \sigma^2(u) \, u_{x} \quad  \text{in } L^2([0,T] \times \R ). \label{eq: L2 conv 0 multi}
\end{align}
To simplify notation, let us temporarily use $f:=\sigma^2(u) \, u_{x}$. Also, for each $ t \in [0,T] $, let 
\begin{align*}
	\widehat{f}(t,z) = \int_\R f(t,x) e^{-2\pi i z x }\, \dd x 
\end{align*}
denote the Fourier transform of $ f(t,\cdot) $. By the Plancherel Theorem,
\begin{align*}
	\int_0^T \! \int_\R \big(f^\varepsilon(t,x) - f(t,x)\big)^2 \, \dd x \, \dd t&= \int_0^T \!  \int_\R \big| \widehat{f}(t,z) \, \widehat{\varphi}_\varepsilon(z) - \widehat{f}(t,z)\big|^2\, \dd z \, \dd t \\
	& = \int_0^T \!  \int_\R  \big| \widehat{f}(t,z) \big|^2 \big(e^{-2\pi^2\varepsilon z^2} - 1\big)^2\, \dd z \, \dd t. 
\end{align*} 

Since $ \big|e^{-2\pi^2\varepsilon z^2} - 1\big| \le 1 $, the claim \eqref{eq: L2 conv 0 multi} follows from the Dominated Convergence Theorem. Moreover,
\begin{align} \label{eq: f_eps L2 conv}
	g^\varepsilon \to g  \quad \text{in } L^2([0,T] \times \R),
\end{align}
where
\begin{align*}
	g^\varepsilon(t, x) &:= \varrho^0_{x}(x - y- b(u^\varepsilon)t - \gamma(u^\varepsilon)  W_t , u^\varepsilon -\eta), \quad 
	g(t,x):= \varrho^0_{x}(x - y- b(u)t  - \gamma(u) W_t , u - \eta).
\end{align*}
Indeed, from Lemma \ref{lem: CDF}(i), Fubini's theorem and Assumption \ref{assum}(a),(b), we see that $g^\varepsilon \to g$  a.e. on $[0,T]  \times \R$. Therefore, \eqref{eq: f_eps L2 conv} follows from the Bounded Convergence Theorem, applied on the common compact support of $ g^\varepsilon $, $ g $. Combining the two $L^2$-convergences \eqref{eq: L2 conv 0 multi} and \eqref{eq: f_eps L2 conv}, we find $ I^\varepsilon_{14} \to I_{14} $. 

	Secondly, we prove that $ I^\varepsilon_{13} \to I_{13} $.   The same argument showing \eqref{eq: L2 conv 0 multi} also shows that $u^\varepsilon_{x}  \to u_{x} $ in $L^2([0,T] \times \R)$. Together with \eqref{eq: L2 conv 0 multi} and the Cauchy-Schwarz inequality, we have $ u_{x}^\varepsilon \left(\sigma^2(u)  u_{x}\right)^\varepsilon \to u_{x} \left(\sigma^2(u)  u_{x}\right) $ in $L^1([0,T] \times \R )$. Therefore, we have the required convergence by Lemma \ref{lem: L1 conv}.

The proof of $ I^\varepsilon_{12} \to I_{12} $ is  the same.  Combining the three convergences $ I^\varepsilon_{14} \to I_{14} $, $ I^\varepsilon_{12} \to I_{12} $ and $ I^\varepsilon_{13} \to I_{13} $, we see that $ \lim_{\varepsilon \downarrow 0} I^\varepsilon_3 = I_3 $ holds. \\ \\
\noindent{\it Group 3: } $ \lim_{\varepsilon \downarrow 0} I^\varepsilon_5 = I_5 $ and $ \lim_{\varepsilon \downarrow 0} I^\varepsilon_8 = I_8 $. \\ \\
We first show that $ \lim_{\varepsilon \downarrow 0} I^\varepsilon_5 = I_5 $.  By the Dambis-Dubins-Schwarz Theorem, it suffices to check that
\begin{align} \label{eq: group3.conv}
\begin{split}
	\lim_{\varepsilon\downarrow 0}\int_s^t  \bigg(\int_\R & \varrho^0\big(x - y- b(u^\varepsilon) r-  \gamma(u^\varepsilon)  W_r , u^\varepsilon -\eta\big) \, G(u)_{x}^\varepsilon \\
&-\varrho^0(x- y - b(u)r - \gamma(u) W_r , u - \eta) \, G(u)_{x} \, \dd x \bigg)^2 \, \dd r  = 0.
\end{split}
\end{align}
Similarly to \eqref{eq: f_eps L2 conv} and \eqref{eq: L2 conv 0 multi}, we can prove that for a.e. $r \in [s,t]$,
\begin{align*}
	&\varrho^0( x - y- b(u^\varepsilon)r - \gamma(u^\varepsilon)  W_r , u^\varepsilon - \eta) 
	\to  \varrho^0(x- y- b(u)r -  \gamma(u) W_r , u- \eta) \hspace{0.7em} \text{in} \hspace{0.7em}  L^2(\R), \\
& G(u)_{x}^\varepsilon \to G(u)_{x} \hspace{0.7em}  \text{in} \hspace{0.7em}  L^2(\R) .
\end{align*}
These two $ L^2 $ convergences imply that the $\dd r$-integrand in \eqref{eq: group3.conv} tends to zero a.e. In conjunction with the Cauchy-Schwarz inequality, $ \Vert G(u)^\varepsilon_x \Vert^2_{L^2(\R)}  \le \Vert G(u)_x \Vert^2_{L^2(\R)} $, $ G(u)_x \in L^2([s,t] \times \R) $, and the Dominated Convergence Theorem, this implies \eqref{eq: group3.conv}.

Similarly, to show that  $ \lim_{\varepsilon \downarrow 0} I^\varepsilon_8 = I_8 $,  it suffices to check that
\begin{align*}
	  \lim_{\varepsilon \downarrow 0} \int_s^t \bigg(\int_{\R} \int_0^1&  \left(\mathbf{1}_{\{ \xi \le u^\varepsilon(r,x)\}} - \mathbf{1}_{\{\xi \le u(r,x)\}}\right) \\
	  &\varrho^0_{x} \left(x- y - b(\xi) r  - \gamma(\xi)  W_r, \xi -\eta \right) \gamma(\xi) \, \dd \xi \, \dd x \bigg)^2 \dd r  = 0. 
\end{align*}
This follows from Lemma \ref{lem: CDF}(i), and two applications of the Bounded Convergence Theorem (recall that $ \varrho^0_x $ is compactly supported, as well as Assumption \ref{assum}(a),(b)). 

\subsection*{Step 6. Proof of \eqref{mainthm: pe}}  As $u(t,\cdot) = F_{\nu_t}(\cdot)$ is a CDF for each $t \in [0,T]$, the requirement $ u \in L^\infty([0,T] \times \R) $ is satisfied.  Also, since $  \nu \in C([0,T]; \P_1(\R))$ by assumption, the integrability and continuity requirements in  \eqref{u.continuity.0}  and \eqref{u.continuity} hold. The condition \eqref{PE.def.1} is fulfilled as $ \sigma $ is bounded and $u_x \in L^2([0,T] \times \R) $ by Step 3.    It remains to show \eqref{eq: pathwise entropy solution result multi} and \eqref{eq: chain rule}.

\subsection*{Step 6.1. Proof of \eqref{eq: pathwise entropy solution result multi}.}  
We  prove \eqref{eq: pathwise entropy solution result multi} with $m \equiv 0$ in this step.  Fix $ (\eta,y) \in \R^2$. As in the previous steps, we  abbreviate $ \varrho (\xi, t,x; \eta ,y)$  in \eqref{eq: varrho multi} by $ \varrho (\xi, t, x)$.  The mollified version of the first term on the LHS of \eqref{eq: pathwise entropy solution result multi} is then
\begin{align*}
	- \int_{\R^{2}} \overline{\chi}(\xi, u^\varepsilon(\cdot ,x)) \,  \varrho( \xi, \cdot,x ) \, \dd \xi \,  \dd x \bigg|_s^t = - \int_{\R} \left( \int_s^t \dd F\left(r, x, u^\varepsilon(r,x), W_r\right)\right) \dd x,
\end{align*} 
where $ F : [0,T] \times \R \times [0,1] \times \R \to \R$ is defined by
\begin{align*}
	F(t, x, u ,w) := \int_0^u \varrho^0(x - y - b(\xi)t -\gamma(\xi)   w, \xi - \eta) \,  \dd \xi. 
\end{align*}
Recalling the dynamics of $ u^\varepsilon $ in \eqref{eq: u eps dym multi},  we have by It\^{o}'s Lemma:
\begin{align*}
	&\dd F\left(r, x, u^\varepsilon(r,x), W_r\right)  \\
	&=-\int_0^{u^\varepsilon} \varrho^0_{x}(x - y- b(\xi)r -\gamma(\xi)  W_r , \xi- \eta) \,  b(\xi) \,  \dd \xi \,  \dd r \\
	& \hspace{0.4cm}- \varrho^0(x - y- b(u^\varepsilon)r - \gamma(u^\varepsilon)  W_r , u^\varepsilon - \eta) 
	\,  \Beta(u)_{x}^\varepsilon \,  \dd r    \\
	& \hspace{0.4cm}+  \varrho^0(x - y - b(u^\varepsilon)r - \gamma(u^\varepsilon)  W_r , u^\varepsilon- \eta)   \,  \Sigma(u)_{xx}^\varepsilon \,   \dd r \\ 
&\hspace{0.4cm}+ \varrho^0(x - y - b(u^\varepsilon)r - \gamma(u^\varepsilon)  W_r , u^\varepsilon - \eta) 
	\,  \Gamma(u)_{xx}^\varepsilon \,  \dd r   \\
	  & \hspace{0.4cm}-   \varrho^0(x  - y- b(u^\varepsilon)r - \gamma(u^\varepsilon) W_r, u^\varepsilon - \eta) \,  G(u)_{x}^\varepsilon\,  \dd W_r \\
	  & \hspace{0.4cm}- \frac{1}{2} \,  \varrho^0_{x }(x - y- b(u^\varepsilon)r - \gamma(u^\varepsilon) W_r , u^\varepsilon - \eta) \,  \big(b'(u^\varepsilon) r + \gamma'(u^\varepsilon) W_r\big) \,   \left(G(u)_{x}^\varepsilon\right)^2 \dd r  \\
	  &\hspace{0.4cm}+ \frac{1}{2} \, \varrho^0_\xi(x  - y- b(u^\varepsilon)r - \gamma(u^\varepsilon) W_r, u^\varepsilon - \eta) \left(G(u)_{x}^\varepsilon\right)^2 \dd r  \\
	  &\hspace{0.4cm}-\int_0^{u^\varepsilon} \varrho^0_{x}(x - y- b(\xi)r  - \gamma(\xi)  W_r , \xi -\eta) \,  \gamma(\xi) \,  \dd \xi \,  \dd W_r\\ 
	 &\hspace{0.4cm}+  \frac{1}{2} \,\int_0^{u^\varepsilon} \varrho^0_{xx}\left(x - y- b(\xi)r - \gamma(\xi)  W_r , \xi - \eta\right) \,  \gamma^2(\xi) \,  \dd \xi \,  \dd r \\
	 &\hspace{0.4cm}+ \varrho^0_{x}(x - y- b(u^\varepsilon)r - \gamma(u^\varepsilon) W_r , u^\varepsilon - \eta) \,  \gamma(u^\varepsilon) \,  G(u)_{x}^\varepsilon \,  \dd r. 
\end{align*}
And so recalling $I^\varepsilon_1, \dots, I^\varepsilon_{10}$ from \eqref{a.to.j.epsilon} and applying Fubini's Theorem and the Stochastic Fubini Theorem (see, e.g., \cite[Theorem 2.2]{veraar2012stochastic}), we get
\begin{align*}
	- \int_{\R^{2}} \overline{\chi}(\xi, u^\varepsilon(\cdot ,x))  \, \varrho( \xi, \cdot, x) \,  \dd \xi \,  \dd x \bigg|_s^t = \sum_{i=1}^{10} I^\varepsilon_i. 
\end{align*}

Together with \eqref{eq: convergence 1} and \eqref{eq: convergence 2} from Step 5 and \eqref{eq: Step 1 result} from Step 4, we have
\begin{align*}
	&- \int_{\R^{{2}}} \overline{\chi}(\xi, u(\cdot ,x))  \, \varrho(\xi, \cdot,x ) \, \dd \xi \, \dd x \bigg|_s^t  +   \frac{1}{2} \int_s^t \int_{\R^{2}} \overline{\chi }\left(\xi, u(r,x)\right)  \sigma^2(\xi ) \, \varrho_{xx}(\xi, r, x) \, \dd \xi \, \dd x \, \dd r \\
	&=-\lim_{\varepsilon \downarrow 0}\int_{\R^{2}} \overline{\chi}(\xi, u^\varepsilon(\cdot ,x))  \varrho(\xi, \cdot, x) \, \dd \xi \, \dd x \bigg|_s^t  + I_{11}\\
	&= \lim_{\varepsilon \downarrow 0} \sum_{i=1}^{10} I^\varepsilon_i + I_{11} = \sum_{i=1}^{10} I_i + I_{11}= I_{12} + I_{13}.
\end{align*}
It remains to show that 
\begin{align} \label{eq:n.equality}
	\int_s^t \int_{\R^{2}}  \partial_{\xi}  \varrho(\xi, r, x) \, n(\dd x, \dd \xi, \dd r)  = I_{12} + I_{13}.
\end{align}
We recall from \eqref{eq: measure n} and \eqref{eq: beta multi} that
\begin{align*}
		n(\dd x, \dd \xi , \dd r) &=  \frac{1}{2} \delta_{u(r,x)}  (\dd \xi ) \, \sigma^2(u(r,x)) \, u_x^2(r,x) \, \dd r \,  \dd x.
\end{align*}
In view of Fubini's Theorem and 
\begin{align*}
	 \partial_{\xi}  \varrho(\xi, r, x)  =& - \varrho_{x}^0 (x - y- b(\xi) r -  \gamma(\xi)W_r , \xi- \eta) \left( b'(\xi)r + \gamma'(\xi) W_r \right) \\
	 & + \varrho_{\xi}^0(x - y- b(\xi)r - \gamma(\xi)W_r , \xi - \eta),
\end{align*}
we see that \eqref{eq:n.equality} holds.  This completes the proof  of \eqref{eq: pathwise entropy solution result multi}. 

\subsection*{Step 6.2. Proof of \eqref{eq: chain rule}}
Let us now prove the ``chain rule" \eqref{eq: chain rule}.  From part \eqref{mainthm: ac} of the theorem, we know that $ u(t,\cdot) $ is continuous for a.e. $t \in [0,T]$.   Fix such a $t \in [0,T]$ and $ (\eta, y) \in \R^2$. 
Recalling the expressions for $ \varrho $ in \eqref{eq: varrho multi} and $ S $ in \eqref{eq: beta multi}, the chain rule becomes
\begin{align*}
	&\int_0^1 \sigma(\xi) \int_{\R} \mathbf{1}_{\{0 \le \xi \le u(t,x)\}} \,  \varrho^0_{x}\big(x - y- b(\xi) t - \gamma(\xi) W_t , \xi - \eta\big) \, \dd x \, \dd \xi \\
	& = - \int_{\R}  \sigma(u(t,x)) \,  u_{x}(t,x) \,  \varrho^0\big(x - y- b(u(t,x))t  -  \gamma(u(t,x)) W_t , u(t,x) - \eta\big) \, \dd x.
\end{align*} 
 For all $ \xi \in (0,1) $,  the inner integral on the LHS is 
\begin{align*}
	&\int_{u^{-1}(t,  \xi)}^\infty \varrho^0_{x}\big(x  -y- b(\xi)t -  \gamma(\xi)  W_t , \xi - \eta\big) \, \dd x  =  - \varrho^0\left(u^{-1}(t, \xi) - y- b(\xi) t - \gamma(\xi)W_t , \xi - \eta\right),
\end{align*}
where $ u^{-1}(t,  \xi) $ is the $\xi$-quantile of $\dd u (t, \cdot)$. Putting this back into the chain rule, it suffices to show that
\begin{align*}
&\int_0^1 \sigma(\xi) \, \varrho^0\left(u^{-1}(t, \xi) - y- b(\xi) t -  \gamma(\xi)W_t , \xi - \eta\right) \dd \xi \\
& =\int_{\R}  \sigma(u(t,x)) \,  u_{x}(t,x) \,  \varrho^0\big(x - y - b(u(t,x))t - \gamma(u(t,x)) W_t , u(t,x) - \eta\big) \, \dd x.
\end{align*}
This can be proven by following the same change-of-variable technique as in \eqref{eq: co-area claim 1 multi}. \hfill \qed

\section{Proof of Theorem \ref{thm: rank}} \label{sec:pf.main.thm}
Using Theorem \ref{thm: main}, we  prove Theorem \ref{thm: rank} in this section.

	{\noindent \bf Step 1. Existence, uniqueness, and tightness.}
	Existence and uniqueness of a weak solution to  \eqref{eq: general rank based system} are consequences of \cite[Theorem 2.1]{bass1987uniqueness}. For the tightness of $\left(\nu^n\right)_{n \in \N}$,  let $\delta  > 0$ and consider two stopping times $0 \le \tau_1 \le \tau_2 \le T$ with $\tau_2 - \tau_1 \le \delta$ a.s. By the Burkholder-Davis-Gundy inequality, there exists a $C < \infty$ such that 
	\begin{align*}
		&\E \big[\big| X^{n,i}_{\tau_2} - X^{n,i}_{\tau_1} \big|\big] \\
	&	\le \E \left[\bigg|\int_{\tau_1}^{\tau_2} b\big(F_{\nu^n_t}\big(X^{n,i}_t\big)\big) \, \dd t \bigg| + \bigg|\int_{\tau_1}^{\tau_2} \sigma \big(F_{\nu^n_t}\big(X^{n,i}_t\big)\big) \, \dd B^{n,i}_t    + \int_{\tau_1}^{\tau_2} \gamma \big(F_{\nu^n_t}\big(X^{n,i}_t\big)\big) \,  \dd W^{n}_t \bigg|  \right] \\
	& 	 \le C \, \E \left[\int_{\tau_1}^{\tau_2} \big| b\big(F_{\nu^n_t}\big(X^{n,i}_t\big)\big) \big| \,\dd t + \sqrt{\int_{\tau_1}^{\tau_2} (\sigma^2 + \gamma^2) \big(F_{\nu^n_t}\big(X^{n,i}_t\big)\big) \, \dd t}  \right] \\
	&	 \le  C \big( \delta \Vert b \Vert_\infty + \sqrt{(\Vert \sigma \Vert_\infty^2 + \Vert \gamma\Vert_\infty^2 ) \delta} \big).
	\end{align*}
Bounding the infimum in the definition of the Wasserstein distance \eqref{eq:def.W1} using the trivial coupling $\frac{1}{n} \sum_{i=1}^n \delta_{(X^{n,i}_{\tau_2}, X^{n,i}_{\tau_1})}$ of $\nu^n_{\tau_2}$ and $\nu^n_{\tau_1}$, we have
	\begin{align*}
		\E \big[\W_1(\nu^n_{\tau_2}, \nu^n_{\tau_1})\big]  \le \frac{1}{n} \sum_{i=1}^{n}  \E \big[\big| X^{n,i}_{\tau_2} - X^{n,i}_{\tau_1} \big|\big] \le C \big( \delta \Vert b \Vert_\infty + \sqrt{(\Vert \sigma \Vert_\infty^2 + \Vert \gamma\Vert_\infty^2 ) \delta} \big),
	\end{align*}
	and so, 
	\begin{align*}
		\lim_{\delta \downarrow 0} \, \limsup_{n \to \infty} \sup_{0 \le \tau_1 \le \tau_2 \le T: \, \tau_2 - \tau_1 < \delta}  \E \big[\W_1(\nu^n_{\tau_2}, \nu^n_{\tau_1})\big]  = 0. 
	\end{align*}
	By Aldous' criterion for tightness \cite[Lemma 23.12, Theorems 23.11, 23.9, 23.8]{kallenberg2021foundations},  we see that the sequence $(\nu^n)_{n \in \N}$  is tight  on $C([0,T];\P_1(\R))  $. \\
	\\{\noindent \bf Step 2. Limit points solve the martingale problem.} 
	From the tightness of $(\nu^n)_{n \in \N}$ in Step 1, we deduce the (joint) tightness of  $(\nu^n, W^{n})_{n \in \N}$ on $C([0,T];\P_1(\R)) \times C([0,T] ; \R)  $. Let $(\nu, W) \in C([0,T] ; \P_1(\R)) \times C([0,T] ; \R)$ be any limit point, supported on some probability space $(\Omega, \F, \PP)$. Note that Assumption \ref{assum:init} ensures that $ \nu_0 = \nu^0 $. We equip  $(\Omega, \F , \PP)$ with the filtration $\FF = \{\F_t\}_{t \in [0,T]}$ generated by $ (\nu, W)$.  In this step, we  show that $ \nu $ induces  a solution to the martingale problem in \eqref{eq: MRT condition}. Fix $k \in \N$ and $f_1, \dots, f_k\in C^\infty_c(\R)$. 
	Let $ \widetilde{f}_i (x) = \int_x^\infty f_i(y)\, \dd y $ for $ i = 1 ,\dots, k $ and note that $ \langle F_{\nu^n_t}, f_i \rangle = \langle \nu^n_t, \widetilde{f}_i  \rangle $.  Together with It\^{o}'s Lemma, we have the dynamics
	\begin{align} \label{eq:ito.dym.1}
	\begin{split}
	\dd  \langle F_{\nu^n_t}, f_i \rangle =& - \Big\langle \nu^n_t , f_i(\cdot) \, b\left(F_{\nu^n_t}(\cdot)\right) +    \frac{1}{2}  f_i' \left(\cdot\right) \big(\sigma^2 + \gamma^2\big)  \left(F_{\nu^n_t}(\cdot)\right)\!\Big\rangle \, \dd t   \\
	& -  \frac{1}{n} \sum_{m=1}^n  f_i \left(X^{n,m}_t\right) \sigma \left(F_{\nu^n_t}\left(X^{n,m}_t\right)\right) \dd B^{n,m}_t \\
	& -  \frac{1}{n} \sum_{m=1}^n  f_i \left(X^{n,m}_t\right) \gamma \left(F_{\nu^n_t}\left(X^{n,m}_t\right)\right) \dd W^{n}_t.
	\end{split}
	\end{align}
	Therefore, the quadratic covariation between $ \langle F_{\nu^n_t}, f_i \rangle$ and $ \langle F_{\nu^n_t}, f_j \rangle$ is given by
	\begin{align} \label{eq:ito.qv}
	\begin{split}
		&\,\,\dd \left\langle \langle F_{\nu^n_t}, f_i \rangle ,\langle F_{\nu^n_t}, f_j \rangle \right\rangle_t\\
	&= \frac{1}{n} \left\langle \nu^n_t, f_i(\cdot) \, f_j(\cdot) \, \sigma^2\left(F_{\nu^n_t}(\cdot)\right)  \right\rangle  \dd t + \left\langle \nu^n_t, f_i(\cdot) \, \gamma \left(F_{\nu^n_t}(\cdot)\right)  \right\rangle \left\langle \nu^n_t, f_j(\cdot) \, \gamma \left(F_{\nu^n_t}(\cdot)\right)  \right\rangle   \dd t.
	\end{split}
	\end{align}
	
	Fix further $\phi \in C^\infty_c(\R^k)$ and $0  \le t \le T$. Using It\^{o}'s Lemma, in conjunction with \eqref{eq:ito.dym.1} and \eqref{eq:ito.qv}, we have
	\begin{align}
	& 	\phi (	\langle F_{\nu^n_t}, \boldsymbol{f} \rangle ) - 	\phi (	\langle F_{\nu^n_0}, \boldsymbol{f} \rangle  ) \label{eq: ito1}\\
	& =-\sum_{i=1}^k \int_0^t  \partial_i \phi (	\langle F_{\nu^n_r}, \boldsymbol{f} \rangle ) \Big\langle \nu^n_r, f_i(\cdot) \, b\left(F_{\nu^n_r}(\cdot)\right) + \frac{1}{2}  f_i' (\cdot)\left(\sigma^2 + \gamma^2\right) \left(F_{\nu^n_r}(\cdot)\right) \!\Big\rangle \,  \dd r  \label{eq: ito2}  \\
	&\hspace{0.5cm}+ \frac{1}{2n} \sum_{i,j=1}^k  \int_0^t  \partial_{ij} \phi (	\langle F_{\nu^n_r}, \boldsymbol{f} \rangle ) \left\langle  \nu^n_r, f_i(\cdot) \, f_j(\cdot) \, \sigma^2\left(F_{\nu^n_r}(\cdot)\right)  \right\rangle  \dd r  \label{eq: ito3}\\
	&\hspace{0.5cm}+ \frac{1}{2} \sum_{i,j=1}^k  \int_0^t \partial_{ij} \phi (	\langle F_{\nu^n_r}, \boldsymbol{f} \rangle ) \left\langle \nu^n_r, f_i(\cdot) \, \gamma \left(F_{\nu^n_r}(\cdot)\right) \right\rangle \left\langle \nu^n_r, f_j(\cdot) \, \gamma \left(F_{\nu^n_r}(\cdot)\right) \right\rangle   \dd r   \label{eq: ito4}\\
	 		&\hspace{0.5cm}-  \frac{1}{n} \sum_{i=1}^k  \sum_{m=1}^n \int_0^t \partial_i \phi (	\langle F_{\nu^n_r}, \boldsymbol{f} \rangle )  f_i \left(X^{n,m}_r\right) \sigma \left(F_{\nu^n_r}\left(X^{n,m}_r\right)\right) \dd B^{n,m}_r  \label{eq: ito5}\\
	&\hspace{0.5cm}- \sum_{i=1}^k \int_0^t \partial_i \phi (	\langle F_{\nu^n_r}, \boldsymbol{f} \rangle ) \left\langle \nu^n_r, f_i (\cdot) \, \gamma \left(F_{\nu^n_r}(\cdot)\right) \right\rangle  \dd W^{n}_r. \label{eq: ito6}
	\end{align} 
	Let us analyze  the subsequential $ n \to \infty $ limits of \eqref{eq: ito1}--\eqref{eq: ito4}. For that purpose, we use the Skorokhod Representation Theorem in the form of \cite[Theorem 3.5.1]{dudley2014uniform} to assume that $(\nu^n, W^n)$ converges a.s. to $(\nu, W)$ on some common probability space. Note that   on this new probability space, each $\nu^n$ admits the representation $\nu^n_r = \frac{1}{n} \sum_{i=1}^n \delta_{X^{n,i}_r}$, where a.s., for a.e. $r \in [0,t]$, the random variables $\{X^{n,i}_r\}_{i=1}^n$ are distinct. Indeed, since on the original probability spaces supporting \eqref{eq: general rank based system}, the diffusion coefficient $\sigma$ is non-degenerate, applying the occupation time formula (see, e.g., \cite[Theorem 3.7.1 and Exercise 3.7.10]{karatzas2012brownian}) to the semimartingale $X^{n,i} - X^{n,j}$ shows that a.s., 
	\begin{align*}
		\int_0^t \mathbf{1}_{\{X^{n,i}_r = X^{n,j}_r \}} \dd r = 0.
	\end{align*}
	\\{\noindent \bf Step 2.1.} 
	For the first term \eqref{eq: ito1}, noting that $\nu^n \to \nu$ a.s., $ \langle F_{\nu^n_t}, f_i \rangle = \langle \nu^n_t, \widetilde{f}_i \rangle $,  $ \langle F_{\nu_t}, f_i \rangle = \langle \nu_t, \widetilde{f}_i \rangle $, and $\phi$ is continuous,
	\begin{align*}
	\phi (\langle F_{\nu^n_t}, \boldsymbol{f} \rangle ) - 	\phi (	\langle F_{\nu^n_0}, \boldsymbol{f} \rangle  )  \to 	\phi (	\langle F_{\nu_t}, \boldsymbol{f} \rangle ) - 	\phi (	\langle F_{\nu_0}, \boldsymbol{f} \rangle  ).
	\end{align*}
	In addition, the third term  \eqref{eq: ito3} converges to zero a.s.  as $n \to \infty$ because the integrand is bounded by Assumption \ref{assum}(c).  
	\\ \\{\noindent \bf Step 2.2.} 
	We show next that the second term \eqref{eq: ito2} and the fourth term \eqref{eq: ito4} converge a.s. to the expected limits: 	
	\begin{align} \label{eq: step 1(b) conv}
		\begin{split}
				& \int_0^t  \partial_i \phi (	\langle F_{\nu^n_r}, \boldsymbol{f} \rangle ) \Big\langle \nu^n_r, f_i(\cdot) \, b\left(F_{\nu^n_r}(\cdot)\right) + \frac{1}{2}  f_i' (\cdot)\left(\sigma^2 + \gamma^2\right) \left(F_{\nu^n_r}(\cdot)\right) \!\Big\rangle \,  \dd r\\
				& \to - \int_0^t     \partial_i\phi \left(	\langle  F_{\nu_r}, \boldsymbol{f}\rangle \right)  \Big(\big\langle \Beta(F_{\nu_r}(\cdot)), f_i'\big\rangle+ \big\langle (\Sigma + \Gamma)(F_{\nu_r}(\cdot)), f_{i}''\big\rangle \Big) \,  \dd r 
		\end{split}
	\end{align}
	and
	\begin{align*}
		& \int_0^t \partial_{ij} \phi (	\langle F_{\nu^n_r}, \boldsymbol{f} \rangle ) \left\langle \nu^n_r, f_i(\cdot) \, \gamma \left(F_{\nu^n_r}(\cdot)\right) \right\rangle \left\langle \nu^n_r, f_j(\cdot) \, \gamma \left(F_{\nu^n_r}(\cdot)\right) \right\rangle   \dd r   \\
		&\to \,    \int_0^t  \partial_{ij} \phi\left(	\langle F_{\nu_r}, \boldsymbol{f} \rangle \right)  \left\langle G(F_{\nu_r}(\cdot)), f_i' \right\rangle \left\langle G(F_{\nu_r}(\cdot)), f_j' \right\rangle  \, \dd r.  
	\end{align*}
	Since both convergences are similar, we focus on \eqref{eq: step 1(b) conv}.  
We claim that  for a.e. $r \in [0,t]$,
	\begin{align}
	& \left\langle \nu^{n}_r, f_i(\cdot) \,  b \left(F_{\nu^{n}_r}(\cdot)\right) \right\rangle  \to -\big\langle \Beta(F_{\nu_r}(\cdot)), f_i'\big\rangle,  \label{eq: f'conv}   \\
	&\left\langle \nu^{n}_r, f_i'(\cdot) \left(\sigma^2 + \gamma^2\right) \left(F_{\nu^{n}_r}(\cdot)\right) \right\rangle  \to -2\big\langle (\Sigma + \Gamma)(F_{\nu_r}(\cdot)), f_{i}''\big\rangle. \label{eq: f''conv}  
	\end{align}
	Since both convergences are similar, we focus on  \eqref{eq: f''conv}.

	For a.e. $ r \in [0,t] $ and distinct $X^{n,1}_r, \dots, X^{n,n}_r$, let us  write 
	\begin{align*}
		\min\{X^{n,1}_r, \dots, X^{n,n}_r\}= X^{n,(1)}_r < X^{n,(2)}_r < \cdots <X^{n,(n)}_r = \max\{X^{n,1}_r, \dots, X^{n,n}_r\}
	\end{align*} for the order statistics. More specifically,  
	\begin{align*}
		X^{n, (\ell)}_r := \min_{1 \le m_1 < \dots < m_\ell \le n} \max\big\{X^{n,m_1}_r, \dots, X^{n, m_\ell}_r\big\}. 
	\end{align*}
	Then, for a.e. $ r \in [0,t] $, 
	\begin{align*}
	& \left\langle \nu^{n}_r, f'_i(\cdot) \left(\sigma^2 + \gamma^2\right) \left(F_{\nu^{n}_r}(\cdot)\right) \right\rangle  = \frac{1}{n} \sum_{m=1}^n f'_i\big(X^{n,m}_r \big)(\sigma^2 + \gamma^2) \big(F_{\nu^{n}_r}(X^{n,m}_r) \big)  \\
	&=\,  \frac{1}{n} \sum_{\ell=1}^n f'_i\big(X^{n,(\ell)}_r\big) (\sigma^2 + \gamma^2) \left(\frac{\ell}{n}\right)   = \frac{1}{n} \int_\R f'_i(y) \,  \dd \sum_{\ell=1}^{n F_{\nu^{n}_r}(y)} (\sigma^2 + \gamma^2) \left(\frac{\ell}{n}\right) \\
	& =- \frac{1}{n} \int_\R \sum_{\ell=1}^{n F_{\nu^{n}_r}(y)} (\sigma^2 + \gamma^2) \left(\frac{\ell}{n}\right) f''_i(y) \,  \dd y.
	\end{align*}
	Therefore, 
	\begin{align}\label{eq: MP drift conv1}
	\begin{split}
		&\left| \left\langle \nu^{n}_r, f'_i(\cdot) \left(\sigma^2 + \gamma^2\right) \left(F_{\nu^{n}_r}(\cdot)\right) \right\rangle  + 2  \int_\R f''_i(y) \, (\Sigma + \Gamma) \!\left(F_{\nu^{n}_r}(y)\right) \dd y  \right| \\
	&= \bigg| \int_\R f''_i(y) \bigg[-\frac{1}{n}  \sum_{\ell=1}^{n F_{\nu^{n}_r}(y)} (\sigma^2 + \gamma^2) \left(\frac{\ell}{n}\right)  +\int_0^{F_{\nu^{n}_r}(y)}(\sigma^2 + \gamma^2)(a) \, \dd a  \bigg] \dd y \bigg| \\
	&\le  \Vert f''_i \Vert_{L^1(\R)} \, \sup_{ y \in \R} \bigg|-\frac{1}{n}  \sum_{\ell=1}^{n F_{\nu^{n}_r}(y)} (\sigma^2 + \gamma^2) \left(\frac{\ell}{n}\right)  +\int_0^{F_{\nu^{n}_r}(y)}(\sigma^2 + \gamma^2)(a) \,  \dd a\bigg| \\
	&= \Vert f''_i \Vert_{L^1(\R)} \, \sup_{q = 1, \dots, n} \bigg|-\frac{1}{n} \sum_{\ell = 1}^q (\sigma^2 + \gamma^2)\left(\frac{\ell}{n} 	\right)+ \int_0^{\frac{q}{n}} (\sigma^2 + \gamma^2)(a) \,  \dd a  \bigg|  \\
	&\le   \Vert f''_i \Vert_{L^1(\R)} \, \sup\left\{\Big|(\sigma^2 + \gamma^2)(a) - (\sigma^2 + \gamma^2)( \widetilde{a})\Big|: a, \widetilde{a} \in [0,1], |a - \widetilde{a}| \le 1/n\right\},
	\end{split}
	\end{align}
	which converges to $ 0 $ as $n \to \infty$ by the uniform continuity of $\sigma$ and $\gamma$.   	On the other hand, since  for all $r \in [0,t]$, $ F_{\nu^{n}_r}  \to F_{\nu_r}  $ a.e., we have by the Dominated Convergence Theorem
	\begin{align}	 \label{eq: MP drift conv2}
	&\lim_{n \to \infty} \int_\R f_i''(y) \, (\Sigma + \Gamma)\! \left(F_{\nu^{n}_r}(y)\right) \dd y =   \int_\R f_i''(y) \,  (\Sigma + \Gamma)\! \left(F_{\nu_r}(y)\right) \dd y.
	\end{align}
	The observations \eqref{eq: MP drift conv1} and \eqref{eq: MP drift conv2} imply  \eqref{eq: f''conv}. Using   \eqref{eq: f'conv}--\eqref{eq: f''conv} in conjunction with 	
 $ \partial_i \phi (\langle F_{\nu^n_r}, \boldsymbol{f} \rangle )  \to   \partial_i \phi (\langle F_{\nu_r}, \boldsymbol{f} \rangle )  $ for all $r \in [0,t]$,    we deduce 
 \begin{align*} 
 		& \partial_i \phi (	\langle F_{\nu^n_r}, \boldsymbol{f} \rangle ) \Big\langle \nu^n_r, f_i(\cdot) \, b\left(F_{\nu^n_r}(\cdot)\right) + \frac{1}{2}  f_i' (\cdot)\left(\sigma^2 + \gamma^2\right) \left(F_{\nu^n_r}(\cdot)\right) \!\Big\rangle \,  \\
 		&\to -     \partial_i\phi \left(	\langle  F_{\nu_r}, \boldsymbol{f}\rangle \right)  \Big(\big\langle \Beta(F_{\nu_r}(\cdot)), f_i'\big\rangle+ \big\langle (\Sigma + \Gamma)(F_{\nu_r}(\cdot)), f_{i}''\big\rangle \Big).
 \end{align*}
	Thus,  \eqref{eq: step 1(b) conv} follows from the Bounded Convergence Theorem. 
	\\\\{\noindent \bf Step 2.3.} 
	We now show that the process $(M_t)_{t \in [0,T]}$ is an $\FF$-martingale, where
	\begin{align*} 
	M_t := &  \,	\phi (	\langle F_{\nu_t}, \boldsymbol{f} \rangle ) - 	\phi (	\langle F_{\nu_0}, \boldsymbol{f} \rangle  )   \\
	\quad & - \sum_{i=1}^k \int_0^t     \partial_i\phi \left(	\langle  F_{\nu_r}, \boldsymbol{f}\rangle \right)  \Big(\big\langle \Beta(F_{\nu_r}(\cdot)), f_i'\big\rangle+ \big\langle (\Sigma + \Gamma)(F_{\nu_r}(\cdot)), f_{i}''\big\rangle \Big) \,  \dd r \\
	\quad & -\frac{1}{2} \sum_{i,j=1}^k \int_0^t  \partial_{ij} \phi\left(	\langle F_{\nu_r}, \boldsymbol{f} \rangle \right)  \left\langle G(F_{\nu_r}(\cdot)), f_i' \right\rangle \left\langle G(F_{\nu_r}(\cdot)), f_j' \right\rangle  \, \dd r.
	\end{align*}
To see this,  first note that $(M^{n}_t)_{t \in [0,T]}$ 	is a  martingale, where 
	\begin{align*}
	M^{n}_t :=  &  -  \frac{1}{n} \sum_{i=1}^k  \sum_{m=1}^n \int_0^t \partial_i \phi(	\langle F_{\nu^n_r}, \boldsymbol{f} \rangle )  f_i \left(X^{n,m}_r\right) \sigma \left(F_{\nu^n_r}\left(X^{n,m}_r\right)\right) \dd B^{n,m}_r   \\
	&  - \sum_{i=1}^k \int_0^t \partial_i \phi (	\langle F_{\nu^n_r}, \boldsymbol{f} \rangle ) \left\langle \nu^n_r, f_i (\cdot) \, \gamma \left(F_{\nu^n_r}(\cdot)\right) \right\rangle  \dd W^{n}_r. 
\end{align*}
 From \eqref{eq: ito1}--\eqref{eq: ito6} and Steps 2.1--2.2,  we see that   for all $t \in [0,T]$,  $M^{n}_t \to M_t  $ a.s. on the new probability space, which implies	 $M^n_t \stackrel{d}{\rightarrow} M_t$ on the original probability spaces.  Similarly, we deduce that the finite-dimensional distributions of $M^{n} $ converge to those of $M$, i.e., $(M^{n}_{t_1}, \dots, M^{n}_{t_\ell}) \stackrel{d}{\rightarrow} (M_{t_1}, \dots, M_{t_\ell})$ for any finite subset $ \{t_1, \dots, t_\ell\} $ of $ [0,T] $. 

	Applying  the Burkholder-Davis-Gundy inequality in the form of  \cite[Exercise 3.3.25]{karatzas2012brownian},  we see that 
			\begin{align*}
	\E \left[\left|M^{n}_{\widetilde{t} } - M^{n}_{t}\right|^4\right]   \le C (\widetilde{t} - t)^2,
	\end{align*}
	where $C$ is a constant depending only on $ \max_{i=1,\dots,k} \Vert f_i \Vert_\infty $, $ \max_{i=1,\dots,k} \Vert \partial_i \phi \Vert_\infty $, $\Vert \sigma \Vert_\infty$ and  $\Vert \gamma \Vert_\infty$.
	Therefore, we conclude from \cite[Problem 2.4.11]{karatzas2012brownian} that $\left(M^{n}\right)_{n \in \N}$ is tight.  As a result, $(M^{n}, \nu^n, W^{n})_{n \in \N}$ is also tight. Hence, $(M^{n}, \nu^n, W^{n}) \stackrel{d}{\to} (M, \nu, W)$ along a subsequence, i.e.,  for  any bounded continuous $\Psi : C([0,T]; \R) \times C([0,T]; \P_1(\R)) \times C([0,T]; \R) \to \R$, we have 
	\begin{align} \label{eq:M.rho.W.weak.conv}
		\E\left[\Psi(M^{n}, \nu^{n }, W^{n})\right] \to \E\left[\Psi(M, \nu, W)\right].
	\end{align}

	 To complete the proof, note that \eqref{eq: ito1}--\eqref{eq: ito6} implies $|M^{n}| \le C$  a.s. for some constant $C < \infty$. Moreover, for any $0 \le s \le t \le T$ and any bounded continuous $\widetilde{\Psi} : C([0,s];\P_1(\R)) \times C([0,s];\R) \to \R$, the function $\Psi : C([0,t];\R) \times C([0,s];\P_1(\R)) \times C([0,s];\R) \to \R$ defined by $\Psi(X,Y,Z) := [(X_t - X_s)\wedge 2C \vee (-2C)] \, \widetilde{\Psi}(Y, Z)$ is also bounded and  continuous.  In conjunction with \eqref{eq:M.rho.W.weak.conv}, we see that 
	\begin{align*} 
	\E \big[(M_t - M_s) \, \widetilde{\Psi}(\nu|_{[0,s]}, W|_{[0,s]})\big] = 
	\lim_{n \to \infty}\E \left[(M^{n}_t - M^{n}_s) \, \widetilde{\Psi}\big(\nu^n|_{[0,s]}, W^{n}|_{[0,s]}\big)\right] = 0,
	\end{align*}
	which shows that $M$ is an $\FF$-martingale.  \\ \\
		\noindent{\bf Step 3. Completing the proof.} 
		 From Step 2 and 
Theorem \ref{thm: main}(i), we see that \eqref{eq: weak sol SPDE} holds. 
 Note that since each $ u(t,\cdot) $ is a CDF, we have $ \Vert u(t, \cdot) \Vert_{BV(\R)} = 1 $.  Therefore, pathwise uniqueness is a consequence of Theorem \ref{thm: main}(iv) and Proposition \ref{prop:GSextension}.  
 
 Finally, for uniqueness in law, consider the following subspace of $  L^\infty([0,T] ; BV(\R))$:
 \begin{align*}
 	\mathcal{S}  := \big\{u \in L^\infty([0,T] ; BV(\R)): u(t, \cdot) = F_{\nu_t}, \nu \in C([0,T]; \P_1(\R))\big\},
 \end{align*}
 equipped with the topology inherited from $  C([0,T]; \P_1(\R)) $. Note that $ \mathcal{S} $ is a Polish space, and is therefore a Borel space on account of \cite[Theorem 1.8]{kallenberg2021foundations}. And so the regular conditional probability for random elements in $ \mathcal{S} $ exists by \cite[Theorem 8.5]{kallenberg2021foundations}.  Finally, $u$ is a random element in $\mathcal{S}$.
  On account of these observations, uniqueness in law for $ u $, and in turn for $ \nu $, follows by a natural extension of the Yamada-Watanabe theorem (see, e.g., \cite[Proposition 5.3.20]{karatzas2012brownian}). \hfill \qed\\

\appendix

\section{Proof of Proposition \ref{prop:GSextension}} \label{sec:proof.gs}
In this section, we  prove Proposition \ref{prop:GSextension}. The proof follows mostly \cite[proof of Theorem 2.3]{gess2017stochastic} on p.~ 2975--2984, which uses some lemmas in their appendix as well. To keep the exposition at a reasonable length, we refer to \cite{gess2017stochastic} for any notations not defined here. There are only a few  main changes that are necessary, so we omit the details that are the same and focus on the differences. They are summarized in the following list. 
\begin{enumerate}[(i)]
	\item Replace every occurrence of $f(\xi)  z_t $ by $ b(\xi)\,t + \gamma(\xi)z_t$. For example, near the bottom of \cite[p.~ 2993]{gess2017stochastic}, the original definition  $\varrho^s_\varepsilon(x, y, \xi, t) := \varrho^s_\varepsilon(x - y + f(\xi)  z_t)$ there is now changed to $\varrho^s_\varepsilon(\xi, t, x; y) := \varrho^s_\varepsilon(x - y + b(\xi) t + \gamma(\xi) z_t)$.  The only exceptions are the equation in the statement of \cite[Lemma A.2]{gess2017stochastic} and the expression in the second equality in the equation display starting with $\tilde{G}(t)$ at the bottom of \cite[p. 2994]{gess2017stochastic}.
	\item Follow the proof of \cite[Theorem 2.3]{gess2017stochastic} up the to the equation after (3.11), i.e., 
	\begin{align}
		G_{\varepsilon, \psi ,\delta}(t) - G_{\varepsilon, \psi ,\delta}(s) \le \int_s^t 	\left(Err^{(1)}_{\varepsilon,\psi,\delta}(r) + Err^{(2)}_{\varepsilon,\psi,\delta}(r) + Err^{(1,2)}_{\varepsilon,\psi,\delta}(r) + Err^{par}_{\varepsilon,\psi,\delta}(r) \right) \dd r& \label{eq:Err12}\\
		+ \int_s^t 	\left(Err^{loc, (1)}_{\varepsilon,\psi,\delta}(r) + Err^{loc,(2)}_{\varepsilon,\psi,\delta}(r) + Err^{loc,(3)}_{\varepsilon,\psi,\delta}(r) + Err^{loc,(4)}_{\varepsilon,\psi,\delta}(r) \right) \dd r&. \label{eq:Errloc1234}
	\end{align}
	Note that the three integrals in $ G_{\varepsilon, \psi, \delta}(t) $, defined in the first display of \cite[p. 2977]{gess2017stochastic}, should be combined into one integral to ensure finiteness by the integrability assumption \eqref{u.continuity.0}. 
	\item For the term in \eqref{eq:Errloc1234}, we first follow the original proof by using \cite[Lemma A.5]{gess2017stochastic} to see that as $ \delta \downarrow 0$, it converges to 
	\begin{align*}
		\int_s^t \left(Err^{loc, (1)}_{\varepsilon,\psi}(r) + Err^{loc,(2)}_{\varepsilon,\psi}(r) + Err^{loc,(3)}_{\varepsilon,\psi}(r) + Err^{loc,(4)}_{\varepsilon,\psi}(r) \right) \dd r. 
	\end{align*}
	Then we begin to diverge from the original proof. We choose  $\psi_R \in C^\infty_c(\R^{2}) $ to be of the form $ \psi_R (\eta, y) = \phi_R(\eta) \, \widetilde{\psi}_R(y) $, where $ \phi_R \in C^\infty_c(\R) $ is such that 
	\begin{align}
	\phi_R(\eta) = \begin{cases}
	1, \quad |\eta| \le R, \\
	0, \quad |\eta| \ge R + 1, 
	\end{cases} |\phi'_R| \le C
	\end{align}
	and $ \widetilde{\psi}_R \in C^\infty_c(\R) $ is such that
	\begin{align*}
	\widetilde{\psi}_R(y) = \begin{cases}
	1, \quad |y| \le R, \\
	0, \quad |y| \ge 2R,
	\end{cases} \quad |\widetilde{\psi}_R'| \le C/R, \quad  |\widetilde{\psi}''_R| \le C/R^{2}
	\end{align*}
	for some $C \ge 1$.   Then we use Lemma \ref{lem:mod.GS.Lem.A.6} below to obtain 	\begin{align}
	\lim_{R \to \infty} \int_s^t \left(Err^{loc, (1)}_{\varepsilon,\psi_R}(r) + Err^{loc,(2)}_{\varepsilon,\psi_R}(r) + Err^{loc,(3)}_{\varepsilon,\psi_R}(r) + Err^{loc,(4)}_{\varepsilon,\psi_R}(r) \right) \dd r = 0. 
	\end{align}
	\item For the term in  \eqref{eq:Err12}, we  follow the original proof by using   \cite[Lemma A.4]{gess2017stochastic} to find
	\begin{align*}
	 \lim_{\delta \downarrow 0} 	\int_s^t Err^{par}_{\varepsilon,\psi,\delta}(r) \, \dd r = 0. 
	\end{align*}
Moreover, by Lemma \ref{lem:mod.GS.A2.A4} below it holds
	\begin{align*}
	\lim_{R \to \infty} \lim_{\delta \downarrow 0}\int_s^t 	\left(Err^{(1)}_{\varepsilon,\psi_R,\delta}(r) + Err^{(2)}_{\varepsilon,\psi_R,\delta}(r) + Err^{(1,2)}_{\varepsilon,\psi_R,\delta}(r)\right) \dd r \le  \frac{C}{\varepsilon}\Vert z^{(1)} - z^{(2)} \Vert_{C([s,t];\R)}.
\end{align*}
\item Combining (iii), (iv), and the fact that $\lim_{R \to \infty} \lim_{\delta \downarrow 0} G_{\varepsilon, \psi_R, \delta}(t) = G_\varepsilon(t)$, we have
\begin{align*}
	G_\varepsilon(t) - G_\varepsilon(s) \lesssim \varepsilon^{-1}  \Vert z^{(1)} - z^{(2)} \Vert_{C([s,t];\R)},
\end{align*} 
which is \cite[equation (3.6)]{gess2017stochastic}. The rest of the proof proceeds as in \cite{gess2017stochastic}. Note that as in  (ii), the three integrals in $ G_\varepsilon(t) $,  defined on \cite[p. 2976]{gess2017stochastic}, should be combined into one integral to  ensure finiteness by the integrability assumption \eqref{u.continuity.0}. 
\end{enumerate}

\begin{lemma} \label{lem:mod.GS.Lem.A.6}
	As $R \to \infty$, 
	\begin{align*}
		\int_s^t Err^{loc,(1)}_{\varepsilon,\psi_R}(r)\,  \dd r \to 0 \quad \text{and} \quad 	\int_s^t \Big(Err^{loc,(3)}_{\varepsilon,\psi_R}(r) +  Err^{loc,(4)}_{\varepsilon,\psi_R}(r)\Big) \dd r \to 0.
	\end{align*} 
\end{lemma} 
\noindent{\bf Proof.}
	To match the notation in \cite{gess2017stochastic}, let us define $a \! : [0,1] \to (0, \infty)$ by $a(r):= \frac{1}{2} \,\sigma^2(r)$. 
	For the first claim, note that $|\partial_{y y}\psi_R(\eta, y)| \le C R^{-2} \mathbf{1}_{[-2R, 2R]}(y)$. Therefore, the integrand in $ Err^{loc,(1)}_{\varepsilon,\psi_R}(r)$, defined on \cite[p. 2997]{gess2017stochastic}, is bounded in absolute value by
	\begin{align*}
	&\big| a(\eta) \,  \partial_{yy}  \psi_R(\eta, y) \, \chi^{(2)}(\eta,r,x') \, \varrho^{s,(2)}_\varepsilon(r,x';\eta,y)  \big| \\ 
	&\le \Vert a \Vert_\infty  \, \frac{C}{R^{2}} \,  \mathbf{1}_{[-2R, 2R]}(y) \,   \varrho^{s,(2)}_\varepsilon(r,x';\eta,y) \, \mathbf{1}_{[-\Vert u^{(2)} \Vert_\infty, \Vert u^{(2)} \Vert_\infty]}(\eta).
	\end{align*}
	Moreover, 
	\begin{align*}
	&\frac{C}{R^{2}}\int_s^t \int_{\R^2} \int_{-\Vert u^{(2)} \Vert_\infty}^{\Vert u^{(2)} \Vert_\infty} \mathbf{1}_{[-2R, 2R]}(y) \, \varrho^{s,(2)}_\varepsilon(r,x';\eta,y) \, \dd \eta \,  \dd x' \, \dd y  \, \dd r \\
	&= \frac{C}{R^{2}} \int_s^t  \int_{\R} \int_{-\Vert u^{(2)} \Vert_\infty}^{\Vert u^{(2)} \Vert_\infty}  \mathbf{1}_{[-2R, 2R]}(y)  \, \dd \eta  \, \dd y \, \dd r\\
	&=   \frac{4C}{R} \int_s^t  \int_{-\Vert u^{(2)} \Vert_\infty}^{\Vert u^{(2)} \Vert_\infty} \dd \eta \, \dd r = \frac{8C \, \Vert u^{(2)} \Vert_\infty \, (t - s)}{R} \to 0,
	\end{align*} 
	which proves the first claim.  
	
	For the second claim, as in \cite[equation (3.9)]{gess2017stochastic}, we can write
	\begin{align*}
	&\varrho^{s,(1)}_\varepsilon(r,x;\eta,y) \, \partial_{x'} \varrho^{s,(2)}_\varepsilon(r,x';\eta,y) + \varrho^{s,(2)}_\varepsilon(r,x';\eta,y) \, \partial_{x} \varrho^{s,(1)}_\varepsilon(r,x;\eta,y) \\
	&= -\varrho^{s,(1)}_\varepsilon(r,x;\eta,y) \,  \partial_{y} \varrho^{s,(2)}_\varepsilon(r,x';\eta,y) - \varrho^{s,(2)}_\varepsilon(r,x';\eta,y) \,  \partial_{y} \varrho^{s,(1)}_\varepsilon(r,x;\eta,y)\\
	&= - \partial_{y} \big(\varrho^{s,(1)}_\varepsilon(r,x;\eta,y) \,  \varrho^{s,(2)}_\varepsilon(r,x';\eta,y)  \big).
	\end{align*}
	Therefore, 
	\begin{align*}
	&\int_s^t \left(Err^{loc,(3)}_{\varepsilon,\psi_R}(r) +  Err^{loc,(4)}_{\varepsilon,\psi_R}(r)\right) \dd r \\
	&=  \int_s^t  \int_{\R^4}  \partial_{y}\psi_R(\eta, y) \, \chi^{(1)}(\eta,r,x) \, \chi^{(2)}(\eta,r,x') \, a(\eta) \, \\
	& \hspace{1.7cm}  \partial_{y} \big(\varrho^{s,(1)}_\varepsilon(r,x;\eta,y) \varrho^{s,(2)}_\varepsilon(r,x';\eta,y)  \big)  \, \dd x \,  \dd x' \, \dd y \, \dd \eta\,  \dd r \\
	&=  - \int_s^t   \int_{\R^4}  \partial_{y y}\psi_R(\eta, y) \, \chi^{(1)}(\eta,r,x) \,  \chi^{(2)}(\eta,r,x') \, a(\eta)\,  \\ &\hspace{2.1cm}  \varrho^{s,(1)}_\varepsilon(r,x;\eta,y) \,  \varrho^{s,(2)}_\varepsilon(r,x';\eta,y)  \, \dd x \, \dd x' \, \dd y \, \dd \eta \, \dd r. 
	\end{align*}
	
	Similarly to the proof of the first claim, the integrand is  bounded in absolute value by 
	\begin{align*}
	&\big| \partial_{yy}\psi_R(\eta, y) \, \chi^{(1)}(\eta,r,x) \,  \chi^{(2)}(\eta,r,x') \,  a(\eta) \,  \varrho^{s,(1)}_\varepsilon(r,x;\eta,y) \,  \varrho^{s,(2)}_\varepsilon(r,x';\eta,y)  \big| \\
	&\le \frac{C}{R^{2}} \mathbf{1}_{[-2R, 2R]}(y) \,  \Vert a \Vert_\infty \,  \varrho^{s,(1)}_\varepsilon(r,x;\eta,y) \,  \varrho^{s,(2)}_\varepsilon(r,x';\eta,y) \,  \mathbf{1}_{[-\Vert u^{(1)} \Vert_\infty, \Vert u^{(1)} \Vert_\infty]}(\eta). 
	\end{align*}
	Moreover, 
	\begin{align*}
	&\frac{C}{R^{2}} 	\int_s^t \int_{\R^3} \int_{-\Vert u^{(1)} \Vert_\infty}^{\Vert u^{(1)} \Vert_\infty} \mathbf{1}_{[-2R, 2R]}(y)  \, \varrho^{s,(1)}_\varepsilon(r,x;\eta,y) \,  \varrho^{s,(2)}_\varepsilon(r,x';\eta,y)  \, \dd \eta  \, \dd x \,  \dd x'  \, \dd y  \,   \dd r \\
	&= \frac{C}{R^{2}} \int_s^t \int_\R \int_{-\Vert u^{(1)} \Vert_\infty}^{\Vert u^{(1)} \Vert_\infty}  \mathbf{1}_{[-2R, 2R]}(y)  \,  \dd \eta \, \dd y \,  \dd r = \frac{8C \, \Vert u^{(1)} \Vert_\infty \, (t - s)}{R} \to 0,
	\end{align*}
	which proves the second claim. \qed

\begin{lemma} \label{lem:mod.GS.A2.A4}There exists a constant $C  < \infty $ such that 
	\begin{align*}
		\lim_{R \to \infty} \lim_{\delta \downarrow 0}\int_s^t 	\left(Err^{(1)}_{\varepsilon,\psi_R,\delta}(r) + Err^{(2)}_{\varepsilon,\psi_R,\delta}(r) + Err^{(1,2)}_{\varepsilon,\psi_R,\delta}(r)\right) \dd r \le  \frac{C}{\varepsilon}\Vert z^{(1)} - z^{(2)} \Vert_{C([s,t];\R)}.
	\end{align*}

\end{lemma}
\noindent{\bf Proof.}
		As the analysis for all three error terms  is similar, we  concentrate on $ Err^{(1)}_{\varepsilon,\psi_R,\delta}(r) $. Recall from \cite[p.~ 2979]{gess2017stochastic} that 
		\begin{align*}
			&\int_s^t  Err^{(1)}_{\varepsilon,\psi_R,\delta}(r) \, \dd r \\
			&=  -\int_{\R^2 \times [s,t] \times \R^4}  \psi_R (\eta, y) \sgn(\xi)\big( \partial_\xi \varrho^{(1)}_{\varepsilon,\delta}(\xi,r,x;\eta,y) \, \varrho^{(2)}_{\varepsilon,\delta}(\xi',r,x';\eta,y) \\ & \hspace{2.7cm}+   \varrho^{(1)}_{\varepsilon,\delta}(\xi,r,x;\eta,y)\, \partial_{\xi'} \varrho^{(2)}_{\varepsilon,\delta}(\xi',r,x';\eta,y)  \big) \, q^{(2)}(\dd x', \dd \xi',\dd r) \,  \dd x \, \dd \xi \,  \dd y \, \dd \eta.  
		\end{align*}
		To begin,  we rewrite  (part of) the integrand as in the first  display of \cite[proof of Lemma A.7]{gess2017stochastic}:
		\begin{align}
			&\partial_\xi \varrho^{(1)}_{\varepsilon,\delta}(\xi,r,x;\eta,y) \, \varrho^{(2)}_{\varepsilon,\delta}(\xi',r,x';\eta,y)  + \varrho^{(1)}_{\varepsilon,\delta}(\xi,r,x;\eta,y) \,\partial_{\xi'} \varrho^{(2)}_{\varepsilon,\delta}(\xi',r,x';\eta,y) \nonumber \\
			&=  \varrho^v_\delta(\xi - \eta) \,  \varrho^v_\delta(\xi' - \eta) \big( \varrho^{s,(1)}_\varepsilon(\xi,r,x;y) \,\partial_{\xi'} \varrho^{s,(2)}_\varepsilon(\xi',r,x';y) + \partial_{\xi}\varrho^{s,(1)}_\varepsilon(\xi,r,x;y)  \,\varrho^{s,(2)}_\varepsilon(\xi',r,x';y) \big) \label{eq:pd.firstterm}\\
			&+ \varrho^{s,(1)}_\varepsilon(\xi,r,x;y) \,\varrho^{s,(2)}_\varepsilon(\xi',r,x';y) \big( \partial_\xi \varrho^v_\delta(\xi - \eta)\, \varrho^v_\delta(\xi' - \eta) + \varrho^v_\delta(\xi - \eta)\, \partial_{\xi'}\varrho^v_\delta(\xi' - \eta)   \big). \label{eq:pd.secondterm}
		\end{align}
		
		Let us study \eqref{eq:pd.secondterm} first. As in the bottom display of \cite[p.~2998]{gess2017stochastic}, integration by parts gives
		\begin{align*}
			&\int_{\R} \psi_R(\eta, y)\, \sgn(\xi)\, \varrho^{s,(1)}_\varepsilon(\xi,r,x;y) \,\varrho^{s,(2)}_\varepsilon(\xi',r,x';y) \\
			& \hspace{0.5cm} \big( \partial_\xi \varrho^v_\delta(\xi - \eta)\, \varrho^v_\delta(\xi' - \eta) + \varrho^v_\delta(\xi - \eta)\, \partial_{\xi'}\varrho^v_\delta(\xi' - \eta)   \big) \, \dd \eta \nonumber\\
			&= \int_\R  \partial_\eta \psi_R(\eta, y)\,  \sgn(\xi)\, \varrho^{s,(1)}_\varepsilon(\xi,r,x;y) \,\varrho^{s,(2)}_\varepsilon(\xi',r,x';y) \,\varrho^v_\delta(\xi - \eta) \,\varrho^v_\delta(\xi' - \eta) \,  \dd \eta.
		\end{align*}
		Therefore, by Fubini's Theorem,
		\begin{align*}
			&\int_{\R^2 \times [s,t] \times \R^4} \psi_R(\eta, y) \, \sgn(\xi)\,  \varrho^{s,(1)}_\varepsilon(\xi,r,x;y) \, \varrho^{s,(2)}_\varepsilon(\xi',r,x';y) \big( \partial_\xi \varrho^v_\delta(\xi - \eta)\,  \varrho^v_\delta(\xi' - \eta) \nonumber \\
			& \hspace{1.9cm} + \varrho^v_\delta(\xi - \eta) \,\partial_{\xi'}\varrho^v_\delta(\xi' - \eta)   \big) q^{(2)}(\dd x', \dd \xi',\dd r) \, \dd x \, \dd \xi \, \dd y \,\dd \eta \nonumber \\
			&= \int_{\R^2 \times [s,t] \times \R^4} \partial_\eta \psi_R(\eta, y) \,  \sgn(\xi) \,  \varrho^{s,(1)}_\varepsilon(\xi,r,x;y) \,  \varrho^{s,(2)}_\varepsilon(\xi',r,x';y) \\
			& \hspace{2.0cm} \varrho^v_\delta(\xi - \eta) \, \varrho^v_\delta(\xi' - \eta) \, q^{(2)}(\dd x', \dd \xi',\dd r) \, \dd x \, \dd \xi\,  \dd y \, \dd \eta \nonumber  \\
			&= \int_{\R^2 \times [s,t] \times \R^3} \partial_\eta \psi_R(\eta, y)  \, \sgn(\xi)  \,\varrho^{s}_\varepsilon(x' - y - b(\xi')r -  \gamma(\xi') \, z^{(2)}_r ) \\
			& \hspace{2.0cm} \varrho^v_\delta(\xi - \eta) \,\varrho^v_\delta(\xi' - \eta) \, q^{(2)}(\dd x', \dd \xi',\dd r) \, \dd \xi  \,\dd y \,\dd \eta. \nonumber
		\end{align*}
	By Fubini's Theorem again, the integral above is bounded in absolute value by
		\begin{align*}
			&\int_{\R^2 \times [s,t]  \times \R^3}   \big|\partial_\eta \psi_R(\eta, y) \big| \, \varrho^{s}_\varepsilon(x'  - y - b(\xi')r - \gamma(\xi') \, z^{(2)}_r) \\
			& \hspace{2.1cm}\varrho^v_\delta(\xi - \eta) \,\varrho^v_\delta(\xi' - \eta) \,q^{(2)}(\dd x', \dd \xi',\dd r) \, \dd \xi \,  \dd y \,\dd \eta \nonumber \\
			& \le  \int_{\R^2 \times [s,t] \times  \R^2}  \big\Vert \partial_\eta \psi_R(\eta, \cdot) \big\Vert_{C(\R; \R)}  \, \varrho^v_\delta(\xi - \eta)\, \varrho^v_\delta(\xi' - \eta) \,q^{(2)}(\dd x', \dd \xi',\dd r) \, \dd \xi \,\dd \eta \\
				& =   \int_{\R^2 \times [s,t] \times \R} \big\Vert \partial_\eta \psi_R(\eta, \cdot) \big\Vert_{C(\R; \R)}  \,  \varrho^v_\delta(\xi' - \eta)\, q^{(2)}(\dd x', \dd \xi',\dd r)  \, \dd \eta.
		\end{align*}
		When  $ \delta \downarrow 0 $, the latter integral converges to 
		\begin{align*}
			\int_{\R^2 \times [s,t]}  \big\Vert \partial_\eta \psi_R(\xi', \cdot) \big\Vert_{C(\R; \R)}  \,  q^{(2)}(\dd x', \dd \xi',\dd r) , 
		\end{align*} 
		which in turn converges to $ 0 $ as $ R \to \infty $ by the Bounded Convergence Theorem (recall that $ q^{(2)}(\dd x', \dd \xi',\dd r) = m^{(2)}(\dd x', \dd \xi',\dd r) + n^{(2)}(\dd x', \dd \xi',\dd r)$, where $ m^{(2)} $ is finite, and $ n^{(2)} $ defined via \eqref{eq: measure n} is also finite by \eqref{PE.def.1}). 
		
		Turning to \eqref{eq:pd.firstterm}, as in \cite[equation (A.1)]{gess2017stochastic}, we write
		\begin{align}
			&\varrho^{s,(1)}_\varepsilon(\xi,r,x;y) \,\partial_{\xi'} \varrho^{s,(2)}_\varepsilon(\xi',r,x';y) + \partial_{\xi}\varrho^{s,(1)}_\varepsilon(\xi,r,x;y)  \,\varrho^{s,(2)}_\varepsilon(\xi',r,x';y) \nonumber \\ 
			&= \,  \varrho^{s,(1)}_\varepsilon(\xi,r,x;y)\,\big( b'(\xi) r + \gamma'(\xi) z^{(2)}_r\big) \,\partial_y \varrho^{s,(2)}_\varepsilon(\xi',r,x';y) \label{eq:pdd.firsttwoterm.1} \\
			&  \hspace{0.5cm}+ \big( b'(\xi) r + \gamma'(\xi) z^{(1)}_r  \big)\, \partial_y \varrho^{s,(1)}_\varepsilon(\xi,r,x;y)\, \varrho^{s,(2)}_\varepsilon(\xi',r,x';y) \label{eq:pdd.firsttwoterm.2} \\
			&   \hspace{0.5cm}- \varrho^{s,(1)}_\varepsilon(\xi,r,x;y) \Big(\big(b'(\xi) - b'(\xi')\big)r + \big(\gamma'(\xi) - \gamma'(\xi')\big) z^{(2)}_r \Big) \partial_y \varrho^{s,(2)}_\varepsilon(\xi',r,x';y) \label{eq:pdd.thirdterm}.
		\end{align}
		For the term in \eqref{eq:pdd.thirdterm}, we have 
		\begin{align*}
			&\int_{\R^2 \times [s,t] \times \R^4} \psi_R(\eta, y) \, \sgn(\xi) \, \varrho^{s,(1)}_\varepsilon(\xi,r,x;y) \Big(\big(b'(\xi) - b'(\xi')\big)r + \big(\gamma'(\xi) - \gamma'(\xi')\big) z^{(2)}_r \Big)  \\
			&\hspace{2cm}   \partial_y \varrho^{s,(2)}_\varepsilon(\xi',r,x';y) \,\varrho^v_\delta(\xi - \eta)  \, \varrho^v_\delta(\xi' - \eta) \,q^{(2)}(\dd x', \dd \xi',\dd r) \,\dd x \,\dd \xi \,\dd y \, \dd \eta \\
			& =  \int_{\R^2 \times [s,t] \times \R^3} \psi_R(\eta, y) \,\sgn(\xi)  \Big(\big(b'(\xi) - b'(\xi')\big)r + \big(\gamma'(\xi) - \gamma'(\xi')\big) z^{(2)}_r \Big)  \\
			& \hspace{2.5cm}    \partial_y \varrho^{s,(2)}_\varepsilon(\xi',r,x';y) \,\varrho^v_\delta(\xi - \eta)   \,\varrho^v_\delta(\xi' - \eta) \,q^{(2)}(\dd x', \dd \xi', \dd r) \,\dd \xi  \,\dd y \,\dd \eta.
		\end{align*}
	By  Fubini's Theorem, the integral is bounded in absolute value by
		\begin{align*}
			&  \int_{\R^2 \times [s,t] \times \R^3} \psi_R(\eta, y)  \Big(r \big|b'(\xi) - b'(\xi') \big| + \Vert z^{(2)} \Vert_{C([s,t]; \R)} \big|\gamma'(\xi) - \gamma'(\xi')\big|  \Big) \big|\partial_y \varrho^{s,(2)}_\varepsilon(\xi',r,x';y) \big| \\
			& \hspace{2.0cm}\varrho^v_\delta(\xi - \eta)  \, \varrho^v_\delta(\xi' - \eta)\, q^{(2)}(\dd x', \dd \xi',\dd r) \, \dd \xi  \,\dd y \,\dd \eta\\
			& \le  \frac{C'}{\varepsilon} \int_{{\R^2 \times [s,t]}} c_{\psi_R, \delta}(\xi') \, q^{(2)} (\dd x', \dd \xi',\dd r),
		\end{align*}
		where $ C' < \infty $ is an absolute constant, and
		\begin{align*}
			 c_{\psi_R, \delta}(\xi'):=   \int_{\R^2}  & \big\Vert \psi_R(\eta, \cdot)\big\Vert_{C(\R; \R)}  \Big(t  \big|b'(\xi) - b'(\xi')\big| + \Vert z^{(2)} \Vert_{C([s,t]; \R)} \, \big|\gamma'(\xi) - \gamma'(\xi')\big| \Big) \\
			 &  \varrho^v_\delta(\xi - \eta)  \, \varrho^v_\delta(\xi' - \eta) \, \dd \eta \,\dd \xi.
		\end{align*}
	Since $ b' $ and $ \gamma' $ are assumed to be continuous, $ c_{\psi_R, \delta}(\xi') $ is bounded uniformly in $ \xi' $ and $ \delta $. Also, $ \lim_{\delta \downarrow 0} c_{\psi_R, \delta} (\xi') = 0$ for each $ \xi' $. Together with $ q^{(2)} = m^{(2)} + n^{(2)}$, \eqref{eq: measure n}, and \eqref{PE.def.1}, this shows that the above integral converges to $ 0 $ as $ \delta \downarrow 0 $ by the Bounded Convergence Theorem. 
			
			The terms in \eqref{eq:pdd.firsttwoterm.1} and \eqref{eq:pdd.firsttwoterm.2} contribute
		\begin{align*}
			&\int_{\R^2 \times [s,t] \times \R^4}  \Big(\varrho^{s,(1)}_\varepsilon(\xi,r,x;y) \, \big(b'(\xi) r + \gamma'(\xi) \, z^{(2)}_r \big)\,\partial_y \varrho^{s,(2)}_\varepsilon(\xi',r,x';y) \\
			& \hspace{2cm}+ \big(b'(\xi) r + \gamma'(\xi) \, z^{(1)}_r \big)\, \partial_y \varrho^{s,(1)}_\varepsilon(\xi,r,x;y) \,\varrho^{s,(2)}_\varepsilon(\xi',r,x';y)\Big) \\
			&\hspace{2.2cm}  \psi_R(\eta, y)  \,\sgn(\xi) \, \varrho^v_\delta(\xi - \eta) \,\varrho^v_\delta(\xi' - \eta)  \,q^{(2)}(\dd x', \dd \xi',\dd r) \,\dd x\, \dd \xi \, \dd y \,\dd \eta.
		\end{align*}
		The following step is similar to the last equation display of \cite[p. 2999]{gess2017stochastic}, which seems to contain a mistake. Integrating by parts in $ y $, the latter integral is equal to 
	\begin{align*}
	&\int_{\R^2 \times [s,t] \times \R^4}   \psi_R(\eta, y) \, \sgn(\xi) \, \varrho^{s,(1)}_\varepsilon(\xi, r, x;y)  \, \gamma'(\xi) \, \partial_y \varrho^{s,(2)}_\varepsilon(\xi', r, x' ; y) \big(z^{(2)}_r - z^{(1)}_r \big)  \\
	&\hspace{2.0cm}  \varrho^v_\delta(\xi - \eta)\, \varrho^v_\delta(\xi' - \eta) \,q^{(2)}(\dd x', \dd \xi',\dd r) \,\dd x \, \dd \xi \, \dd y \, \dd \eta\\
	&- \int_{\R^2 \times [s,t] \times \R^4}   \partial_y\psi_R(\eta, y) \,  \sgn(\xi)\, \varrho^{s,(1)}_\varepsilon(\xi, r, x; y)  \big(b'(\xi) r + \gamma'(\xi)\,z^{(1)}_r  \big)  \varrho^{s,(2)}_\varepsilon(\xi', r, x'; y) \\
	& \hspace{2.4cm}  \varrho^v_\delta(\xi - \eta)\, \varrho^v_\delta(\xi' - \eta)\, q^{(2)}(\dd x', \dd \xi' ,\dd r) \,\dd x \, \dd \xi \, \dd y \, \dd \eta\\
&=\int_{\R^2 \times [s,t] \times \R^3}   \psi_R(\eta, y) \, \sgn(\xi) \,  \gamma'(\xi) \, \partial_y \varrho^{s,(2)}_\varepsilon(\xi', r, x' ; y) \big(z^{(2)}_r - z^{(1)}_r \big)  \\
&\hspace{2.4cm}  \varrho^v_\delta(\xi - \eta)\, \varrho^v_\delta(\xi' - \eta) \,q^{(2)}(\dd x', \dd \xi',\dd r)  \, \dd \xi \,\dd y \, \dd \eta\\
&\hspace{0.5cm}- \int_{\R^2 \times [s,t] \times \R^3}   \partial_y\psi_R(\eta, y) \,  \sgn(\xi)\,  \big(b'(\xi) r + \gamma'(\xi)\,z^{(1)}_r  \big)  \varrho^{s,(2)}_\varepsilon(\xi', r, x'; y) \\
& \hspace{2.9cm}  \varrho^v_\delta(\xi - \eta)\, \varrho^v_\delta(\xi' - \eta)\, q^{(2)}(\dd x', \dd \xi' ,\dd r)   \, \dd \xi \,\dd y \, \dd \eta.\\
	\end{align*}
	The integrals are bounded in absolute value by 
\begin{align}
	&\Vert z^{(1)} - z^{(2)} \Vert_{C([s,t];\R)}  \int_{\R^2 \times [s,t] \times \R^3} \psi_R(\eta, y)  \big|\gamma'(\xi) \big| \big| \partial_y \varrho^{s,(2)}_\varepsilon(\xi', r, x'; y) \big| \nonumber  \,  \varrho^v_\delta(\xi - \eta)\, \varrho^v_\delta(\xi' - \eta)\\
	& \hspace{5.4cm}   q^{(2)}(\dd x', \dd \xi',\dd r)  \, \dd \xi \,\dd y \, \dd \eta \nonumber\\
	&+   \big( \Vert b' \Vert_\infty \,t+ \Vert\gamma'\Vert_\infty \, \Vert z^{(1)}\Vert_{C([s,t];\R)}  \big)\int_{\R^2 \times [s,t]  \times \R^3} \big|\partial_y\psi_R(\eta, y) \big|   \varrho^{s,(2)}_\varepsilon(\xi', r, x';y) \, \varrho^v_\delta(\xi - \eta)\, \nonumber \\
	& \hspace{8cm} \varrho^v_\delta(\xi' - \eta)\,q^{(2)}(\dd x', \dd \xi',\dd r)  \, \dd \xi \,\dd y \, \dd \eta\nonumber\\
&	\le  \frac{C'}{\varepsilon}\Vert z^{(1)} - z^{(2)} \Vert_{C([s,t];\R)} \, \Vert \gamma' \Vert_\infty \int_{\R^2 \times [s,t] \times \R^2}   \varrho^v_\delta(\xi - \eta) \, \varrho^v_\delta(\xi' - \eta)\, q^{(2)}(\dd x', \dd \xi',\dd r)  \, \dd \xi \, \dd \eta \nonumber\\
	& \hspace{0.4cm}+  \frac{C}{R}  \big( \Vert b' \Vert_\infty \,t+ \Vert\gamma'\Vert_\infty \, \Vert z^{(1)}\Vert_{C([s,t];\R)}  \big)  \int_{\R^2 \times [s,t] \times  \R^2}  \varrho^v_\delta(\xi - \eta) \, \varrho^v_\delta(\xi' - \eta)\, q^{(2)}(\dd x', \dd \xi',\dd r) \, \dd \xi \, \dd \eta\nonumber\\
	&\le  \frac{C_1}{\varepsilon}  \Vert z^{(1)} - z^{(2)} \Vert_{C([s,t];\R)} + \frac{C_2}{R}  \big( \Vert b' \Vert_\infty \,t+ \Vert\gamma'\Vert_\infty \, \Vert z^{(1)}\Vert_{C([s,t];\R)}  \big) , \label{eq: boundB.62}
\end{align}
where 
\begin{align*}
	C_1 &:=   C'  \, q^{(2)}(\R^2 \times [s,t]) \, \Vert \gamma'\Vert_\infty , \quad
	C_2 :=    C \,q^{(2)}(\R^2 \times [s,t] ) .
\end{align*}
Note that both $C_1$ and $ C_2 $ are finite by $ q^{(2)} = m^{(2)} + n^{(2)}$, \eqref{eq: measure n}, and \eqref{PE.def.1}.  Therefore, letting $R \to \infty$, we see that only the first term in \eqref{eq: boundB.62} remains. 

All in all, we see that 
	\begin{align*}
	\lim_{R \to \infty} \lim_{\delta \downarrow 0} \bigg|\int_s^t 	 Err^{(1)}_{\varepsilon,\psi_R,\delta}(r)  \dd r \bigg|\le \frac{C_1 }{\varepsilon}  \Vert z^{(1)} - z^{(2)} \Vert_{C([s,t];\R)},
\end{align*}
as desired. \qed

\section{Some auxiliary results}\label{sec:aux}
Recall the definition of $\varphi_{\varepsilon}$ given in \eqref{eq: varphi}.
\begin{lemma} \label{lem: CDF}
	Let $ F : \R \to [0,1]$ be a CDF and $ f \in L^1(\R)$. Then:
	\begin{enumerate}[(i)]
		\item $ F^\varepsilon(x) \to F(x)$ except at countably many $ x $.
		\item $ f^\varepsilon \to f $ in $ L^1(\R) $. 
	\end{enumerate}
\end{lemma}
\noindent{\bf Proof.}
		Since $ \dd F^\varepsilon $ converges weakly to $ \dd F $, claim {\itshape (i)} follows from \cite[Proposition 3.2.18]{dembonotes}.
		For claim {\itshape (ii)},  note that
		\begin{align*}
			\int_\R \big| f^\varepsilon(x) - f(x) \big|\,  \dd x \le \int_\R \varphi_\varepsilon(y) \int_\R |f(x -y) - f(x) | \, \dd x \, \dd y.
		\end{align*}
		The function $ y \mapsto \int_\R |f(x -y) - f(x) | \, \dd x $ is bounded and has value $ 0 $ at $ 0 $. It is also continuous by the Kolmogorov-Riesz Theorem \cite[Theorem 5(iii)]{hanche2010kolmogorov}. These observations show claim {\itshape (ii)}.  \qed

\begin{lemma} \label{lem: L1 conv}
	Let $f$ and $\{f_k\}_{k \in \N}$ be uniformly bounded functions on $[0,T] \times \R $ such that $f_k \to f$ a.e. Also, let $ g $ and $\{g_k\}_{k \in \N} $ be $L^1([0,T] \times \R)$ functions such that $g_k \to g$ in $L^1([0,T] \times \R)$.  Then, 
	\begin{align*}
 \lim_{k \to \infty}	\int_{[0,T] \times \R} f_k \, g_k   \,\dd t \, \dd x= \int_{[0,T] \times \R} f \, g  \, \dd t \, \dd x.
	\end{align*}
\end{lemma}
\noindent{\bf Proof.}
 We have
		\begin{align*}
		&\bigg|\int_{[0,T] \times \R}  f_k \,g_k  \, \dd t \,\dd x- \int_{[0,T] \times \R} f  g\,\dd t \, \dd x \bigg| \\
		&\le \bigg|\int_{[0,T] \times \R} f_k \, (g_k - g) \, \dd t \, \dd x \bigg| + \bigg|\int_{[0,T] \times \R} (f_k - f) \, g  \, \dd t \,\dd x\bigg|  \\
		&\le \Vert f_k \Vert_\infty \Vert g_k - g\Vert_{L^1([0,T] \times \R)} + \int_{[0,T] \times \R}  |f_k - f| \, |g|  \,\dd t \, \dd x.
		\end{align*}
		The first term converges to zero by assumption, and the second term tends to zero by the Dominated Convergence Theorem. \qed

\begin{lemma}\label{lem:moll.com}
	For each $ \varepsilon > 0 $ and $ x \in \R $, there exists $ (g_m)_{m \in \N} \subset C^\infty_c(\R)$ such that $ g^{(k)}_m \to \varphi^{(k)}_{\varepsilon}(\cdot - x) $ in $ L^1(\R) $ for $ k = 0,1,2 $, where the superscript denotes the $ k $-th order derivative.
\end{lemma}
\noindent{\bf Proof.}
		Let $\psi_m \in C^\infty_c(\R)$ be such that	\begin{align*}
			\psi_m(y) = \begin{cases}
				1, \quad |y| \le m, \\
				0, \quad |y| \ge m+1,
			\end{cases} \quad |\psi_m^{(k)}| \le C^{(k)} \quad \text{for } k =0, 1,2
		\end{align*}
		and some $C^{(k)} < \infty$, $ k = 0,1,2 $. Set $g_m(y) = \varphi_{\varepsilon}(y - x) \, \psi_m(y)$.  It is readily checked that the sequence  $ (g_m )_{m \in \N}$ satisfies the claim.  \hfill \qed

\bibliographystyle{amsplain} 
\bibliography{biblio}
\end{document}